\documentclass[11pt]{amsart}

\usepackage{mathrsfs}
\usepackage{amsmath, amsfonts, amsthm, amssymb}
\usepackage[all,cmtip,line,arc]{xy}
\usepackage{array}
\usepackage{enumerate}
\usepackage{xcolor}
\usepackage{srcltx} 

\newtheorem{thm}{Theorem}[section]
\newtheorem{cor}[thm]{Corollary}
\newtheorem{prop}[thm]{Proposition}
\newtheorem{lem}[thm]{Lemma}

\theoremstyle{definition}
\newtheorem{defn}[thm]{Definition}

\theoremstyle{remark}
\newtheorem{rem}[thm]{Remark}

\newcommand \eps{\varepsilon}

\makeatletter
\let\c@equation\c@thm
\makeatother
\numberwithin{equation}{section}

\title{Interior H\"older estimate for the linearized complex Monge-Ampere equation}
\author{Yulun Xu }
\date{December 2022}

\begin{document}

\maketitle
\begin{abstract}
    Let $w_0$ be a bounded,  $C^3$,  strictly plurisubharmonic function defined on $B_1\subset \mathbb{C}^n$.  Then $w_0$ has a neighborhood in $L^{\infty}(B_1)$. Suppose that we have a function $\phi$ in this neighborhood with $1-\eps \le MA(u)\le 1+\eps$ and there exists a function $u$ solving the linearized complex Monge-Ampere equation: $det(\phi_{k\Bar{l}})\phi^{i\Bar{j}}u_{i\Bar{j}}=0$. Then one has an estimate on $|u|_{C^{\alpha}(B_{\frac{1}{2}})}$ for some $\alpha>0$ depending on $n$, as long as $\epsilon$ is small depending on $n$. This partially generalizes Caffarelli's estimate for linearized real Monge-Ampere equation to the complex version.
\end{abstract}

\section{introduction}

Monge-Ampere equations are second-order partial differential equations whose leading term is the determinant of the Hessian of a real unknown function.  The Hessian is required to be positive or at least nonnegative,  so that the equations are elliptic or degenerate elliptic.  Monge-Ampere equations can be divided into real or complex,  depending on whether one is considering real Hessian or complex Hessian.  In the real case,  the Hessian is $\phi_{ij}$,  so that the positivity of the Hessian is a convexity condition.  In the complex case,  the Hessian is $\phi_{i\bar{j}}$,  and its positivity is a plurisubharmonicity condition.

Let $\phi$ be a convex solution to a real Monge-Ampere equation:
\begin{equation}\label{real ma}
    det D^2 \phi =g.
\end{equation}

\begin{defn}
Let $E\subset \mathbb{C}^n$ be a set and $x_0\in E$.  We will sometimes denote $E$ to be $E(x_0)$ to indicate it is a ``pointed set".  Let $c>0$,  we define:
\begin{equation*}
cE(x_0)=\{x_0+c(y-x_0):y\in E(x_0)\}.
\end{equation*}
Namely $cE(x_0)$ is the image of the dilation map centered at $x_0$ by a factor $c$.
\end{defn}

\begin{defn}
    Let $\mu$ be the Monge-Ampere measure of $\phi$, (In the case of $\phi\in C^2$, $\mu(A)=\int_A det(D^2 \phi)$ for any set $A$) We say $\mu$ satisfies the doubling property if there exist constants $C>0$ and $0<\alpha <1$ such that:
    \begin{equation*}
        \mu(S_t(x)) \le C\mu(\alpha S_t (x)),
    \end{equation*}
    for any section
    \begin{equation*}
        S_t(x) =\{y\in R^n: \phi(y)<l(y)+t\},
    \end{equation*}
    where $l$ is a supporting hyperplane of $\phi$ at $x$. 
\end{defn}
Note that if we have that $\lambda<g < \Lambda$ for some positive constants $\lambda$ and $\Lambda$, then the doubling property holds. 

Next we consider the following linearized Monge-Ampere equation:
\begin{equation*}
    L_{\phi}= det(D^2 \phi) \phi^{ij}u_{ij}=f.
\end{equation*}
When we take first derivatives of $\phi$, we can see that $\phi_j=D_j \phi$, $j=1,...,n$ satisfy the linearized Monge-Ampere equation:
\begin{equation*}
    L_{\phi}(\phi_j)=g_j.
\end{equation*}
Since $\phi$ is convex, the linearized Monge-Ampere equation is elliptic. However, the linearized Monge-Ampere equation is not uniformly elliptic unless we have the estimate for the second derivatives of $\phi$. The standard H\"older estimates for the solutions to linear second order elliptic equations usually require the uniform ellipticity. However, Caffarelli prove the H\"older estimate for the solutions to the linearized Monge-Ampere equations under a weak condition on $g$ which doesn't imply the uniform ellipticity of the linearized Monge-Ampere equation, see\cite{CG2}:
\begin{thm}
    Assume that the Monge-Ampere measure $\mu$ satisfies the doubling property. Let $u$ be a nonnegative solution to the equation:
    \begin{equation*}
        L_{\phi}u=0
    \end{equation*}
    in a section $S_R(x_0)$. Then there exist constants $C_0>0$ and $\alpha>0$ depending on $n$ and $|u|_{\infty}$ and the constants in the doubling property such that:
    \begin{equation*}
        ||u||_{C^{\alpha}(S_{\frac{R}{2}}(x_0))}\le C_0.
    \end{equation*}
\end{thm}

The boundary Harnack inequality for the linearized real Monge-Ampere equation is derived in \cite{L}. There are estimates for the high order derivatives of the solutions to the linearized real Monge-Ampere equation. When $g$ is continuous, the $C^{1,\alpha}$ estimate is derived in \cite{GCT2} and the $W^{2,p}$ estimate is derived in \cite{GCT1}. The boundary H\"older gradient estimates is derived in \cite{LO}. If $g$ is not continuous but belongs to some VMO-type space, the interior $W^{2,p}$ estimate is derived in \cite{H} while the global $W^{2,p}$ estimate is derived in \cite{LZ1}. The $C^{1,\alpha}$ estimate is derived in \cite{LZ2}.

There are some applications of the Real linearized Monge-Ampere equation to the complex geometry. It can be used to prove the interior regularity of the Calabi flow on a toric surface, see \cite{CHL}. It can also be applied to the extremal metrics on toric surfaces, see \cite{Z}. However, as far as I am concerned, the theory of the real linearized Monge-Ampere equation can only be applied to the toric case where a complex Monge-Ampere equation can be reduced to a real Monge-Ampere equation. Besides, the complex linearized Monge-Ampere equation appears in the complex geometry such as the study of the csck problem \cite{CC}. So a natural question is: how to adapt the method for the real linearized Monge-Ampere equation to the complex linerized Monge-Ampere equation directly? Thanks to \cite{CX}, we can give a partial answer to this question:  

\begin{thm}\label{harnack all scales}
    Let $\Omega\subset \mathbb{C}^n$ be a bounded domain with $B_{1-\gamma_0}\subset \Omega\subset B_{1+\gamma_0}$. Let $\phi\in C^2(\Omega)\cap PSH(\Omega)\cap C(\bar{\Omega})$ be such that $1-\eps \le \det \phi_{i\bar{j}}\le 1+\eps$ in $\Omega$ and $\phi=0$ on $\partial \Omega$. Suppose that $\gamma_0$ and $\epsilon$ are small constants depending on $n$. Let $S_t(z_0)$ be defined in \cite{CX}. Then there exist constants $\beta>1$, $\mu_0$ and $C_0$ depending on $n$. Suppose that $u\in C^2(S_{t}(z_0))$ is a nonnegative solution to $L_{\phi}u=0$ on $S_{t}(z_0)$ with $t\le \frac{\mu_0^5}{C_0}$ and $z_0\in B_{\frac{1}{2}}(0)$. Then we have that:
    \begin{equation*}
        \sup_{S_{t}(z_0)}u \le \beta \inf_{S_{t}(z_0)}u.
    \end{equation*}
\end{thm}

\begin{cor}\label{main theorem baby}
Let $\Omega\subset \mathbb{C}^n$ be a bounded domain with $B_{1-\gamma_0}\subset \Omega\subset B_{1+\gamma_0}$. Let $\phi\in C^2(\Omega)\cap PSH(\Omega)\cap C(\bar{\Omega})$ be such that $1-\eps \le \det \phi_{i\bar{j}}\le 1+\eps$ in $\Omega$ and $\phi=0$ on $\partial \Omega$. Suppose that $\gamma_0$ and $\epsilon$ are small constants depending on $n$. Let $S_t(z_0)$ be defined in \cite{CX}. Suppose that $u\in C^2(\Omega)$ is a solution to $L_{\phi}u=0$ on $\Omega$. Then we have that:
\begin{equation*}
    ||u||_{C^{\alpha}(B_{\frac{1}{2}})}\le C,
\end{equation*}
Here $\alpha>0$ is a constant depending on $n$ and  $C$ is a constant depending on $n$ and $|u|_{L^{\infty}(\Omega)}$.
\end{cor}

More generally, we have that:
\begin{thm}\label{main theorem}
Let $w_0$ be a smooth function in the unit ball such that for some $C_0>1$:
\begin{equation*}
      \frac{1}{C_0}I\le (w_0)_{z_i\bar{z}_j}\le C_0I,\,\,\,|D^3w_0|\le C_0\text{ in $B_1$.}
\end{equation*}
Then there exists $\delta_0>0$ small enough,  depending only on $C_0$ and $n$,  such that for all $\phi\in C^2(B_1)\cap PSH(B_1)\cap C(B_1)$ with $|\phi-w_0|\le \delta_0$ on $B_1$,  solving $1-\eps\le MA(\phi)\le 1+\eps$, and for any solution $u\in C^2(B_1)$ solving 
\begin{equation*}
    L_{\phi}u =det(\phi_{k\Bar{l}})\phi^{i\Bar{j}}u_{i\Bar{j}}=0,
\end{equation*}
we have that:
\begin{equation*}
    ||u||_{C^{\alpha}(B_{\frac{1}{2}})}\le C.
\end{equation*}
Here $\alpha>0$ is a constant depending on $n$. Here $C$ is a constant depending on $C_0$, $|u|_{L^{\infty}(B_1)}$ and $n$. $\eps$ is small enough depending only on $n$.
\end{thm}

In the above,  $MA(\phi)$ is the complex Monge-Ampere operator defined for continuous plurisubharmonic functions,  in the Bedford-Taylor sense (see \cite{BT}),  so that $MA(\phi)=\det \phi_{i\bar{j}}$ when $\phi\in C^2$. From now on, we use $L_{\phi}$ for the complex linearized Monge-Ampere equation.

For the manifold setting, we have the following Corollary:
\begin{cor}\label{c1.2}
 Let $(M,\omega_0)$ be a compact K\"ahler manifold.  
Let $\phi\in C^2(M) \cap PSH(M,\omega_0)$ be the solution to:
\begin{equation*}
(\omega_0+\sqrt{-1}\partial\bar{\partial}\phi)^n=f\omega_0^n,\,\,\,\omega_0+\sqrt{-1}\partial\bar{\partial}\phi>0,
\end{equation*}
where $|f-1|<\eps$ and $\int_Mf\omega_0^n=\int_M\omega_0^n$.
Let $u\in C^2(M)$ be the solution to the equation:
\begin{equation*}
    \Delta_{\phi}u=g_{\phi}^{i\Bar{j}} u_{i\Bar{j}}=n-tr_{g_{\phi}}g
\end{equation*}
Suppose that $\epsilon$ is small enough depending on $n$, $\omega_0$. Then we have that
\begin{equation*}
    ||u||_{C^{\alpha}}\le C,
\end{equation*}
Here $\alpha>0$ is a constant depending on $n$ and $C$ is a constant depending on $n$, $\omega_0$ and $||u||_{L^{\infty}}$.
\end{cor}

In the section 3, we reduce the Theorem \ref{main theorem} to the Corollary \ref{main theorem baby}. Then we prove the Corollary \ref{main theorem baby} starting by proving a version of Calderon-Zygmund decomposition in the section 4(Theorem \ref{cz}). Then we prove that the level sets of solutions have uniform critical density in the section 5(Theorem \ref{critical density}). Then we prove that solutions that are large on a section are uniformly large on a bigger section in the section 6(Theorem \ref{larger set inf estimate}). In the section 7,  we first prove the power decay of the distribution function of solutions and then prove the Harnack inequality (Theorem \ref{initial harnack}) and the H\"older estimate of the solutions (Corollary \ref{main theorem baby}). In the section 8, we prove some corollaries of the main theorem.

\section{preliminary}
We want to show that the equations in the main theorem is invariant under affine transformations. Let $z$ be the original coordinate. For any affine transformation $T$ and any positive constant $\lambda$, we can define a new coordinate $w$ by $z=\sqrt{\lambda}Tw$. Let $h$ be a degree two pluriharmonic polynomial. Then we can normalize $\phi$ and $u$ by:
\begin{equation*}
    \begin{split}
        & \widetilde{\phi}(w)=\frac{\phi(\sqrt{\lambda}T w)}{\lambda |det(T)|^{\frac{2}{n}}} +h \\
        & \widetilde{u}(w)=u(\sqrt{\lambda}T w).
    \end{split}
\end{equation*}
Then by calculation, we have that:
\begin{equation*}
    L_{\widetilde{\phi}}\widetilde{u}(w)=\lambda |det(T)|^{\frac{2}{n}} L_{\phi}u(z).
\end{equation*}
So if $L_{\phi}u=0$, we can get that $L_{\widetilde{\phi}}\widetilde{u}=0$
Recall that We denote the complex Monge-Ampere measure as $\mu=MA(\phi)$. We denote the Lebesgue measure as $m$.

\section{Reduction of Theorem \ref{main theorem} to Corollary \ref{main theorem baby}}
We first need the following lemmas from \cite{CX}:
\begin{lem}\label{l2.1}
Let $w_0$ be as stated in Theorem \ref{main theorem}.  Denote $a_{x_0,ij}=(w_0)_{i\Bar{j}}(x_0)$ and $h_{x_0}=Re(\Sigma_i 2(w_0)_i z_i) +Re(\Sigma_{i,j}(w_0)_{ij}z_i z_j)$. Namely we assume that $w_0\in C^3(B_1)$,  and $\frac{1}{C_0}I\le (w_0)_{i\bar{j}}\le C_0I$,  $|D^3w_0|\le C_0$ on $B_{0.99}$.  
Let $\delta\ge 0$ and $\phi_0$ be a function on $B_1$ with $|\phi_0-w_0|\le \delta$ on $B_{0.95}$.  Then there exists $C_1>0$ large enough and $\mu_0>0$ small enough depending only on $C_0$,  such that for all $\mu$ with $4C_1\delta\le \mu\le \mu_0$,  we have: 
\begin{equation*}
(1-C_1\gamma)E_{\mu}(x_0)\subset \{z\in B_{\frac{1}{2C^2_0}}(x_0): (\phi_0-h_{x_0})(z)\le \phi_0(x_0)+\mu\}\subset (1+C_1\gamma)E_{\mu}(x_0).
\end{equation*}
Moreover,  $(\phi_0-h_{x_0})(z)=\phi_0(x_0)+\mu$ on $\partial \{z\in B_{\frac{1}{2C^2_0}}(x_0): (\phi_0-h_{x_0})(z)\le \phi_0(x_0)+\mu\}$.
\end{lem}
Here $\gamma=\frac{\delta}{\mu}+\mu^{\frac{1}{2}}$ and $E_{\mu}(x_0)=\{z\in\mathbb{C}^n:\sum_{i,j=1}^na_{x_0,ij}(z-x_0)_i\overline{(z-x_0)_j}\le \mu\}$.

Let $\phi$ and $w_0$ be as stated in Theorem \ref{main theorem}.
Let $\mu>0$ and $x_0\in B_{0.8}$.  Let $T_{\mu,x_0}$ be a $\mathbb{C}$-linear transformation such that $T_{\mu,x_0}0=0$ and $x_0 +T_{\mu,x_0}(B_{\sqrt{\mu}}(0))=E_{\mu}(x_0)$.  Define 
\begin{equation}\label{2.3N}
\phi_{\mu,x_0}(\zeta)=\frac{1}{\mu|\det T_{\mu,x_0}|^{\frac{2}{n}}}(\phi-h_{x_0}-\mu)(x_0+ T_{\mu,x_0}(\sqrt{\mu}\zeta)).
\end{equation}
Since $E_{\mu}(x_0)$ is defined in terms of $a_{x_0,ij}$,  with $\frac{1}{C_0}\le a_{x_0,ij}=(w_0)_{i\Bar{j}}(x_0)\le C_0I$,  it is easy to see that:
\begin{equation*}
||T_{\mu,x_0}||\le C_2,\,\,\,||T_{\mu,x_0}^{-1}||\le C_2,\,\,\,\frac{1}{C_2}\le |\det T_{\mu,x_0}|^2\le C_2.
\end{equation*}
Here $C_2$ is a large enough constant depending only on $C_0$ and $n$.
Define $\Omega_{\mu}=T_{\mu,x_0}^{-1}(\{z\in B_{\frac{1}{2C_0^2}}:(\phi-h_{x_0})(z)\le \phi(x_0)+\mu\}-x_0)$.
Then by straightforward calculation and Lemma \ref{l2.1},  we can see the following:
\begin{lem}\label{L3.2N}
There is $\mu_0>0$ small enough depending only on $C_0$ such that for all $4C_1\delta_0\le \mu\le \mu_0$ (with $C_1>0$ being the constant given by Lemma \ref{l2.1}),  we have
\begin{enumerate}
\item $B_{1-C_1\gamma}\subset \Omega_{\mu}\subset B_{1+C_1\gamma}$,  with $\gamma=\frac{\delta_0}{\mu}+\mu^{\frac{1}{2}}$.
\item $\det(\phi_{\mu,x_0})_{\zeta_i\bar{\zeta}_j}=f(x_0+T_{\mu,x_0}(\sqrt{\mu}\zeta))$ in $\Omega_{\mu}$,  $u_{\mu,x_0}=0$ on $\partial\Omega_{\mu}$.
\end{enumerate}
\end{lem}
The renormalized function $u_{\mu,x_0}$ fits in the assumptions for Theorem \ref{main theorem baby} after suitably choosing the parameters,  and Theorem \ref{main theorem} follows as a direct consequence:
\begin{cor}\label{c2.3}
Theorem \ref{main theorem} holds,  if we assume Corollary \ref{main theorem baby}.
\end{cor}

\begin{proof}
We wish to apply Corollary \ref{main theorem baby} to each $u_{\mu,x_0}$.  In order to do so,  we just need:
\begin{equation*}
C_1\gamma=C_1(\frac{\delta_0}{\mu}+\mu^{\frac{1}{2}})\le \gamma_0(n),\,\,\,|f(x_0+T_{\mu,x_0}(\sqrt{\mu}\zeta))-1|\le \eps(n).
\end{equation*}
Here $\gamma_0(n)$ and $\eps(n)$ are the constants given by Theorem \ref{main theorem baby}.

So we could just take $\mu$ so that $2C_1\mu^{\frac{1}{2}}\le \frac{1}{2}\gamma_0(n)$ and also $\mu\le \mu_0$ (given by Lemma \ref{L3.2N}).  With this $\mu$,  we can take $\delta_0$ so that $C_1\frac{\delta_0}{\mu}\le \frac{1}{2}\gamma_0(n)$ and also that $4C_1\delta_0\le \mu$.  We fix this choice from now on. 

Since we assumed that Corollary \ref{main theorem baby} holds,  we conclude that:
\begin{equation*}
||u_{\mu,x_0}||_{C^{\alpha}}(B_{\frac{1}{2}})\le C,
\end{equation*}
where $C$ is a constant depending only on $n$ and $|u|_{L^{\infty}(B_1)}$.  Then using (\ref{2.3N}) we may go back to $u$ and obtain that 
\begin{equation*}
||u||_{C^{\alpha}(E_{\frac{1}{2}\mu}(x_0))}\le C',
\end{equation*}
for any $x_0\in B_{\frac{1}{2}}(0)$ where $C'$ is a constant depending only on $n$, $C_0$ (defined in the statement of the Theorem \ref{main theorem}) and $|u|_{L^{\infty}(B_1)}$. Since there exists a constant $C_1$ depending on $C_0$ such that 
\begin{equation*}
   B_{\frac{\sqrt{\mu}}{C}} \subset E_{\frac{1}{2}\mu}(x_0) \subset B_{C\sqrt{\mu}}(x_0).
\end{equation*}
We can use en elementary covering argument to get that:
\begin{equation*}
    ||u||_{C^{\alpha}(B_{\frac{1}{2}}(0))}\le C''.
\end{equation*}
$C''$ is a constant depending only on $n$, $C_0$ and $|u|_{L^{\infty}(B_1)}$.
\end{proof}

\section{Calderon-Zygmund decomposition}
From now on we focus on the Corollary \ref{main theorem baby}. In this section we want to prove the following theorem which is a version of Calderon-Zygmund decomposition using the sections $S_{t}(x)$ defined in \cite{CX}.
\begin{thm}\label{cz}
    Let $0< \sigma<1$ and $0< \delta<1$ be given. Let $\mu_0$ (This is the constant we use to define sections in the $W^{2,p}$ paper) and $\epsilon$ (This is the same constant in the Theorem \ref{main theorem baby}) be small depending on $\sigma$, $\delta$ and $n$. Let $A$ be a bounded subset of $B_{\frac{1}{2}}(0)$. Suppose that for a.e. $x\in A$, 
    \begin{enumerate}
    \item $\lim_{t \rightarrow 0}\frac{\mu(S_t(x)\cap A)}{\mu(S_t(x))}=1.$
    \item $\mu(S_t(x)\cap A)\le  \delta \mu(S_t(x))$ for any $\mu_0^4< t \le \mu_0^3$.
    \end{enumerate}
     Then for such $x$, we can define $t_x=\sup\{t\le \mu_0^4: \mu(S_t(x)\cap A)\ge \delta \mu(S_t(x))\}$. 
    
    Then there exist a countable family of sections $\{S_k=S_{t_k}(x_k)\}$, where $x_k\in A$ and $t_k\le \mu_0^3$, with the following properties:

    (a) $(1-C\sigma^{\frac{1}{2}})\delta\le \frac{\mu(S_k\cap A)}{\mu(S_k)}\le \delta$. 

    (b) For a.e. $x\in A$, $x\in \cup_k S_k$.

    (c) $\mu(A)\le \delta_0 \mu(\cup_1^{\infty}S_k)$, where $\delta_0=\delta_0(\delta)<1$.
\end{thm}
\begin{rem}
    As can be seen from the proof, $t_k$ may not be equal to $t_{x_k}$.
\end{rem}

First we need the following proposition from the \cite{CX} paper.
\begin{prop}\label{section similar to a ball 2}
Let $\Omega$ and $u$ be as stated in Corollary \ref{main theorem baby},  with $\gamma_0$ small enough depending only on  $n$.  Let $0<\sigma<1$ be given.  Then there exists $\eps>0$ depending only on $\sigma$ and $n$,  such that if $|f-1|\le \eps$,  the following hold:
\begin{enumerate}
\item There exists $\mu_0>0$ small enough depending only on $n$ and $\sigma$,  such that for all $x_0\in B_{0,8}$ and all $\mu\le \mu_0$,  there exists a degree 2 pluriharmonic polynomial $h_{\mu,x_0}(z)$ with $h_{\mu,x_0}(x_0)=0$,  such that
\begin{equation*}
\begin{split}
&(1-0.1\sigma)E_{\mu}(x_0)\subset S_{\mu}(x_0):=\{z\in \Omega: (u-h_{\mu,x_0})(z)\le u(x_0)+\mu\}\\
&\subset (1+0.1\sigma)E_{\mu}(x_0).
\end{split}
\end{equation*}
In the above,  $E_{\mu}(x_0)=\{z\in\mathbb{C}^n:\sum_{i,j=1}^na_{\mu,x_0,ij}(z-x_0)_i\overline{(z-x_0)_j}\le \mu\}$,  with $a_{\mu,x_0,ij}$ being positive Hermitian and $\det a_{\mu,x_0,ij}=1$.
\item There is a function $c(\sigma):\sigma\in (0,1)\rightarrow \mathbb{R}_{>0}$,  such that for any $x_0\in B_{0.8}$ and any $0<\mu_1\le \mu_2\le \frac{\mu_0}{1+c(\sigma)}$,  one has $S_{\mu_1}(x_0)\subset S_{(1+c(\sigma))\mu_2}(x_0)$.  Moreover,  $0<c(\sigma)\le C_{2,n}\sigma^{\frac{1}{2}}$ for some dimensional constant $C_{2,n}$.  
\item There is a dimensional constant $C_{3,n}>0$ such that for all $0<\mu\le \mu_0$ and any $x_0\in B_{0.8}$,  there exists a $\mathbb{C}$-linear transformation $T_{\mu,x_0}$,  such that $|\det T_{\mu,x_0}|=1$, $T_{\mu,x_0}0=0$, $x_0 + T_{\mu,x_0}(B_{\sqrt{\mu}}(0))=E_{\mu}(x_0)$.  
Moreover,  for any $0<\mu_1<\mu_2\le \mu_0$ and any $x_0\in B_{0.8}$:
\begin{equation*}
||T_{\mu_1,x_0}\circ T_{\mu_2,x_0}^{-1}||\le C_{3,n}(\frac{\mu_2}{\mu_1})^{\frac{C_{3,n}\sigma^{\frac{1}{2}}}{-\log(0.1\sigma)}},\,\,\,||T_{\mu_2,x_0}\circ T_{\mu_1,x_0}^{-1}||\le C_{3,n}(\frac{\mu_2}{\mu_1})^{\frac{C_{3,n}\sigma^{\frac{1}{2}}}{-\log(0.1\sigma)}}.
\end{equation*}
\end{enumerate}
\end{prop}

The following conclusion is proved in the "induction hypothesis" part of the \cite{CX} paper.
\begin{lem}\label{tx0k estimate}
    For any $\mu_0^{k+1}< t \le \mu_0^k$, we can write 
    \begin{equation*}
        T_{\mu,x_0}=T_{x_0,k}=\widetilde{T}_{x_0,1} \circ \widetilde{T}_{x_0,2}\circ ... \circ \widetilde{T}_{x_0,k},
    \end{equation*}
    where $T_{\mu,x_0}$ is used in the statement of the Proposition \ref{section similar to a ball 2}. we have that 
    \begin{equation*}
        \begin{split}
            &|\widetilde{T}_{x_0,1}|\le C,\,\,\, |\widetilde{T}_{x_0,1}^{-1}|\le C \\
            & |\widetilde{T}_{x_0,k}-I|\le C\sigma^{\frac{1}{2}} \text{ for } k\ge 2.
        \end{split}
    \end{equation*}
\end{lem}

We need the following engulfing property of the sections which is proved in the \cite{CX} paper.
\begin{prop}\label{engulfing property}
Assume that $x_1,\,x_2\in B_{0.8}$,  $0<\mu_1,\,\mu_2\le \mu_0$ and $\mu_1\le 4\mu_2$.  Let $\sigma>0$ be small enough (depending only on dimension).  Assume also that $S_{\mu_1}(x_1)\cap S_{\mu_2}(x_2)\neq \emptyset$,  then $S_{\mu_1}(x_1)\subset 10S_{\mu_2}(x_2)$.
\end{prop}

There is another version of the engulfing property:

\begin{lem}\label{engulfing property theta version}
    Let $\sigma$ be small and $\epsilon$ be small. There exists a constant $\theta>0$ such that if $S_t(z)$ is a section with $y\in S_t(z)$, then $S_t(z)\subset S_{\theta t}(y)$ for $t\le \frac{\mu_0}{\theta}$.
    \end{lem}
\begin{proof}
    Using the Proposition \ref{section similar to a ball 2} we can get that 
    \begin{equation}\label{compare two sections with ellipsoids}
    \begin{split}
         &(1-\sigma)\sqrt{t}E_1\subset S_{t}(y) \subset (1+\sigma)\sqrt{t}E_1 \\
        & (1-\sigma)\sqrt{t_2}E_2\subset S_{t_2}(y) \subset (1+\sigma)\sqrt{t_2}E_2,
    \end{split}
    \end{equation}
    where $E_1$ and $E_2$ are ellipsoids centered at $z$. Let $k_1$ and $k_2$ be integers such that:
    \begin{equation*}
        \mu_0^{k_1+1}<t \le \mu_0^{k_1}, \,\,\, \mu_0^{k_2+1}<t_2\le \mu_0^{k_2}.
    \end{equation*}
    If $t$ and $t_2$ are in the same generation or adjacent generations, i.e. $k_1-1\le k_2\le k_1+1$. So by the Lemma \ref{tx0k estimate}, we have that 
    \begin{equation*}
        |T_{y,k_2} \circ T_{y,k_1}^{-1} -I|\le C \sigma^{\frac{1}{2}},\,\,\,\, |T_{y,k_1} \circ T_{y,k_2}^{-1} -I|\le C \sigma^{\frac{1}{2}}.
    \end{equation*}
    Then we can define $T=T_{y,k_1} \circ T_{y,k_2}^{-1}$ such that $|T-I|\le C\sigma^{\frac{1}{2}}$ and $T E_2 =E_1$. So we have that:
    \begin{equation*}
    \begin{split}
        &10S_{t}(y)\subset 10\sqrt{t}(1+\sigma)E_1 \subset 10\sqrt{t}(1+\sigma)(1+C\sigma^{\frac{1}{2}})E_2 \\
        &=(1-\sigma)\sqrt{100(1+\sigma)^2(1+C\sigma^{\frac{1}{2}})^2\frac{1}{(1-\sigma)^2}t}E_2 \subset S_{100(1+\sigma)^2(1+C\sigma^{\frac{1}{2}})^2\frac{1}{(1-\sigma)^2}t}(y)
    \end{split}
    \end{equation*}
    Then we use the Proposition \ref{engulfing property} to get that:
    \begin{equation*}
        S_{t}(z) \subset 10 S_{t}(y).
    \end{equation*}
    So we have that:
    \begin{equation*}
        S_{t}(z) \subset S_{100(1+\sigma)^2(1+C\sigma^{\frac{1}{2}})^2\frac{1}{(1-\sigma)^2}t}(y).
    \end{equation*}
    We can assume that $\sigma\le \frac{1}{2}$ and let $\mu_0$ be small such that $t$ and $100(1+\sigma)^2(1+C\sigma^{\frac{1}{2}})^2\frac{1}{(1-\sigma)^2}t$ are in the same generation or adjacent generations. In conclusion we can take $\theta=100(1+\sigma)^2(1+C\sigma^{\frac{1}{2}})^2\frac{1}{(1-\sigma)^2}$.  Then we finish the proof the lemma.
\end{proof}

We also need the following lemma from the \cite{CX} estimating the shape of the sections in the original coordinate. In particular, we can get an estimate of the diameter of the sections.
\begin{lem}\label{diam limit}
Suppose that $\mu_0$ is small and $\epsilon$ is small. We have that:
\begin{equation*}
 B_{\frac{1}{C}\mu^{\frac{1}{2}+log_{\mu_0}(1-C\sigma^{\frac{1}{2}})}}(x_0) \subset S_{\mu}(x_0) \subset B_{C \mu^{\frac{1}{2}+log_{\mu_0} (1+C\sigma^{\frac{1}{2}})}}(x_0),
\end{equation*}
for any $\mu\le \mu_0^2$.
\end{lem}
\begin{proof}
We can assume that $\mu \le \mu_0^2$. For any $\mu$, there exists an    integer $k$ such that 
\begin{equation*}
    \mu_0^{k+1}< \mu \le \mu_0^k. 
\end{equation*}
As in the proof of the Proposition \ref{section similar to a ball 2} and the Lemma \ref{tx0k estimate} we have that 
\begin{equation}\label{3.21}
\begin{split}
        &x_0 +(1-\sigma)\widetilde{T}_{x_0,1} \circ \widetilde{T}_{x_0,2}\circ ... \circ \widetilde{T}_{x_0,k} (B_{\sqrt{\mu}}(0)) \subset S_{\mu}(x_0)  \\
        &\subset x_0+ (1+\sigma) \widetilde{T}_{x_0,1} \circ \widetilde{T}_{x_0,2}\circ ... \circ \widetilde{T}_{x_0,k} (B_{\sqrt{\mu}}(0)).
\end{split}
\end{equation}
Using the Lemma \ref{tx0k estimate}, we have that
 \begin{equation*}
        \begin{split}
            &|\widetilde{T}_{x_0,1}|\le C,\,\,\, |\widetilde{T}_{x_0,1}^{-1}|\le C \\
            & |\widetilde{T}_{x_0,k}-I|\le C\sigma^{\frac{1}{2}} \text{ for } k\ge 2.
        \end{split}
    \end{equation*}
So the Formula \ref{3.21} becomes:
\begin{equation*}
\begin{split}
& B_{\frac{1}{C}\mu^{\frac{1}{2}+log_{\mu_0}(1-C\sigma^{\frac{1}{2}})}}(x_0)\subset \frac{1}{C}(1-\sigma) (1-C\sigma^\frac{1}{2})^{k-2}B_{\sqrt{\mu}}(x_0) \subset S_{\mu}(x_0)\\
& \subset C(1+\sigma)(1+C\sigma^{\frac{1}{2}})^{k-2} B_{\sqrt{\mu}}(x_0)\subset B_{C \mu^{\frac{1}{2}+log_{\mu_0} (1+C\sigma^{\frac{1}{2}})}}(x_0).
\end{split}
\end{equation*}
This concludes the proof of the lemma.
\end{proof}

The following lemma is a characterization of sections. The lemma follows directly from the construction of the sections.
\begin{lem}\label{characterization of sections}
    For any $x_0\in B_{0.8}$ and $\mu\le \mu_0$, there exists a degree 2 pluriharmonic polynomial $h_{x_0,\mu}$ such that $h_{x_0,\mu}(x_0)=0$ and
    \begin{equation*}
        S_{\mu}(x_0)=\{x: \phi(x)\le h_{x_0,\mu}(x)+\phi(x_0)+\mu\}.
    \end{equation*}
\end{lem}

The following is the key lemma in the \cite{CX}. It basically means that two sections of the same height of two functions which differ by a plurisubharmonic function are close to each other.
\begin{lem}\label{key lemma 2}
Let $\phi$ be a function defined on an open set $U\subset \mathbb{C}^n$ and let $h(z)$ be pluriharmonic function on $U$.  Let $0<\mu<1$ and $\sigma>0$ be such that:
\begin{equation*}
\begin{split}
&(1-\gamma)\widetilde{E}_{\sqrt{\sigma}}\subset \{\phi\le \sigma\}\subset (1+\gamma)\widetilde{E}_{\sqrt{\sigma}}\subset U,\\
&(1-\gamma)E_{\sqrt{\sigma}}\subset\{\phi \le h+\sigma\}\subset (1+\gamma)E_{\sqrt{\sigma}}\subset U.
\end{split}
\end{equation*}
In the above,  $\widetilde{E}_{\sqrt{\sigma}}=\{z\in\mathbb{C}^n:\sum_{i,j=1}^n \widetilde{a}_{i\Bar{j}}z_i \Bar{z}_j\le \sigma\}$  and $E_{\sqrt{\sigma}}=\{z\in\mathbb{C}^n:\sum_{i,j=1}^n a_{i\Bar{j}}z_i \Bar{z}_j\le \sigma\}$.  Then there exists $c_1(\gamma)$ which is universal (depending only on dimension and you can explicitly calculate) and $c_1(\gamma)\rightarrow 0$ as $\gamma\rightarrow 0$ such that \begin{equation*}
    (1-c_1(\gamma)) \widetilde{E}_{\sqrt{\sigma}} \subset E_{\sqrt{\sigma}} \subset (1+c_1(\gamma)) \widetilde{E}_{\sqrt{\sigma}}
\end{equation*}
\end{lem}

Now we want to prove the following lemma which implies that if two sections have nonempty intersection, then they are comparable to each other:
\begin{lem}\label{advanced engulfing property}
    Let $S_{t_0}(x_0)$ and $S_t(x)$ be two sections such that $t\le t_0\le \frac{\mu_0}{2}$ and 
    \begin{equation*}
        S_{t_0}(x_0) \cap S_t(x) \neq \emptyset.
    \end{equation*}
    Let $T_{t_0,x_0}$ be the affine transformation defined in the $W^{2,p}$ paper that normalize $S_{t_0}(x_0)$.
    Then we have that
    \begin{equation*}
          \frac{1}{C}B_{(\frac{t}{t_0})^{\frac{1}{2}+\epsilon_1}}(\frac{1}{\sqrt{t_0}}T_{t_0,x_0}^{-1} x) \subset \frac{1}{\sqrt{t_0}}T_{t_0,x_0}^{-1} S_{t}(x) \subset C B_{(\frac{t}{t_0})^{\frac{1}{2}-\epsilon_1}}(\frac{1}{\sqrt{t_0}}T_{t_0,x_0}^{-1} x).
    \end{equation*}
    $\epsilon_1$ is a positive constant that can be made arbitrarily small if we let $\mu_0$ be small enough. The constant $C$ depends on $n$. 
\end{lem}
\begin{proof}
     We can use the part (3) of the Proposition \ref{section similar to a ball 2} to get that 
    \begin{equation*}
        S_{t_0}(x_0) \subset S_{2t_0}(x_0),\,\,\, S_{t}(x) \subset S_{2t_0}(x).
    \end{equation*}
    Since $S_{t_0}(x_0)\cap S_t(x)\neq \emptyset$, we have that $S_{2t_0}(x_0)\cap S_{2t_0}(x)\neq \emptyset$. By the Proposition \ref{engulfing property} we have that:
    \begin{equation*}
        S_{2t_0}(x_0) \subset 10 S_{2t_0}(x),\,\,\, S_{2t_0}(x)\subset 10 S_{2t_0}(x_0).
    \end{equation*}
    So we have that:
    \begin{equation*}
        \begin{split}
            &(1-\sigma)(x_0+T_{2t_0,x_0} B_{\sqrt{2t_0}}(0))\subset (1-\sigma)E_{2t_0}(x_0) \subset S_{2t_0}(x_0)\\
            &\subset 10 S_{2t_0}(x) \subset 10(1+\sigma)E_{2t_0}(x)=10(2+\sigma)(x+T_{2t_0,x}B_{\sqrt{2t_0}}(0)). \\
            &(1-\sigma)(x+T_{2t_0,x} B_{\sqrt{2t_0}}(0))\subset (1-\sigma)E_{2t_0}(x) \subset S_{2t_0}(x)\\
            &\subset 10 S_{2t_0}(x_0) \subset 10(1+\sigma)E_{2t_0}(x_0)=10(2+\sigma)(x_0+T_{2t_0,x_0}B_{\sqrt{2t_0}}(0)).
        \end{split}
    \end{equation*}
    So $T_{2t_0,x_0}$ and $T_{2t_0,x}$ are bounded from each other. I.e. 
    \begin{equation}\label{T2t0x0 T2t0x comparable}
        |T_{2t_0,x_0}^{-1} \circ T_{2t_0,x}|\le C, |T_{2t_0,x}^{-1} \circ T_{2t_0,x_0}|\le C.
    \end{equation}
    By the Lemma \ref{tx0k estimate}, $T_{2t_0,x}$ and $T_{t_0,x}$ differ by a linear transformation which is $C\sigma^{\frac{1}{2}}-$close to $Id$. So we have that $T_{2t_0,x}$ and $T_{t_0,x}$ are bounded from each other. Similarly, $T_{2t_0,x_0}$ and $T_{t_0,x_0}$ are bounded from each other. In conclusion, $T_{t_0,x}$ and $T_{t_0,x_0}$ are bounded from each other.
    We consider the following two cases:
    
    (1) $t\le 2 t_0 {\mu_0}^2$. In a coordinate $w$ where $S_{2t_0}(x)$ is close to a ball, we can define $\widetilde{S}_t(x)$ just like how we define $S_t(x)$. In the same coordinate, using the Lemma \ref{diam limit}, we can get that:
\begin{equation*}
 B_{\frac{1}{C}(\frac{t}{2t_0})^{\frac{1}{2}+log_{\mu_0}(1-C\sigma^{\frac{1}{2}})}}(0) \subset \widetilde{S}_{t}(x) \subset B_{C (\frac{t}{2t_0})^{\frac{1}{2}+log_{\mu_0} (1+C\sigma^{\frac{1}{2}})}}(0).
\end{equation*}
Here we use $\frac{t}{2t_0}$ because in the coordinate $w$, $S_{2t_0}(x)$ is close to the unit ball. So the height of the section $\widetilde{S}_{t}(x)$ is scaled to $\frac{t}{2t_0}$ accordingly. By the Lemma \ref{key lemma 2}, $S_t(x)$ is comparable to $\widetilde{S}_t(x)$, i.e. 
\begin{equation*}
    \frac{1}{C}S_t(x) \subset \widetilde{S}_t(x) \subset C S_t(x).
\end{equation*}
So we have that
\begin{equation*}
 B_{\frac{1}{C}(\frac{t}{2t_0})^{\frac{1}{2}+log_{\mu_0}(1-C\sigma^{\frac{1}{2}})}}(0) \subset S_{t}(x) \subset B_{C (\frac{t}{2t_0})^{\frac{1}{2}+log_{\mu_0} (1+C\sigma^{\frac{1}{2}})}}(0).
\end{equation*}
    If we go back to the original coordinate, we have that:
      \begin{equation*}
      B_{\frac{1}{C}(\frac{t}{2t_0})^{\frac{1}{2}+log_{\mu_0}(1-C\sigma^{\frac{1}{2}})}}(0) \subset \frac{1}{\sqrt{2t_0}}T_{2t_0,x}^{-1}(S_t(x)-x)  \subset B_{C (\frac{t}{2t_0})^{\frac{1}{2}+log_{\mu_0} (1+C\sigma^{\frac{1}{2}})}}(0).  
    \end{equation*}
 We already prove that $T_{2t_0,x}$ and $T_{2t_0,x_0}$ are bounded from each other.
    \begin{equation*}
    \begin{split}
        &B_{\frac{1}{C}(\frac{t}{2t_0})^{\frac{1}{2}+log_{\mu_0}(1-C\sigma^{\frac{1}{2}})}}(0) \subset \frac{1}{\sqrt{2t_0}}T_{2t_0,x_0}^{-1}(S_t(x)-x)  \\
        &\subset B_{C (\frac{t}{2t_0})^{\frac{1}{2}+log_{\mu_0} (1+C\sigma^{\frac{1}{2}})}}(0).  
    \end{split}
    \end{equation*}
    So we can take $\epsilon_1=\max{\{-log_{\mu_0} (1+C\sigma^{\frac{1}{2}}), log_{\mu_0}(1-C\sigma^{\frac{1}{2}})\}}$ to get that:
     \begin{equation*}
      B_{\frac{1}{C}(\frac{t}{2t_0})^{\frac{1}{2}+\epsilon_1}}(0) \subset \frac{1}{\sqrt{2t_0}}T_{2t_0,x_0}^{-1}(S_t(x)-x)  \subset B_{C (\frac{t}{2t_0})^{\frac{1}{2}-\epsilon_1}}(0).  
    \end{equation*}
    Recall that $T_{2t_0,x_0}$ and $T_{t_0,x_0}$ are bounded from each other. So we have that
        \begin{equation*}
          \frac{1}{C}B_{(\frac{t}{t_0})^{\frac{1}{2}+\epsilon_1}}(\frac{1}{\sqrt{t_0}}T_{t_0,x_0}^{-1} x) \subset \frac{1}{\sqrt{t_0}}T_{t_0,x_0}^{-1} S_{t}(x) \subset C B_{(\frac{t}{t_0})^{\frac{1}{2}-\epsilon_1}}(\frac{1}{\sqrt{t_0}}T_{t_0,x_0}^{-1} x).
    \end{equation*}
    
    (2) $t_0 \mu_0^2< t \le t_0$. By the Lemma \ref{section similar to a ball 2} we have that:
    \begin{equation*}
        (1-\gamma)B_{\sqrt{t}}(0) \subset T_{t,x}^{-1}(S_{t}(x)-x) \subset (1+\gamma) B_{\sqrt{t}}(0).
    \end{equation*}
    By the Lemma \ref{tx0k estimate}, we have that: $T_{2t_0,x_0}$ and $T_{t,x_0}$ are bounded from each other. $T_{2t_0,x}$ and $T_{t,x}$ are bounded from each other. Combining these facts and the Inequalities \ref{T2t0x0 T2t0x comparable} we have that $T_{t,x}$ and $T_{t_0,x_0}$ are bounded from each other. So we have that:
    \begin{equation*}
        B_{\frac{1}{C}\sqrt{\frac{t}{t_0}}}(\frac{1}{\sqrt{t_0}}T_{t_0,x_0}^{-1}x) \subset \frac{1}{\sqrt{t_0}}T_{t_0,x_0}^{-1} S_t(x) \subset B_{C\sqrt{\frac{t}{t_0}}}(\frac{1}{\sqrt{t_0}}T_{t_0,x_0}^{-1}x).
    \end{equation*}
\end{proof}

Consider the following Dirichlet problem on a domain $\Omega$.
\begin{equation}\label{3.1}
\begin{split}
   & \det((v_0)_{i\bar{j}})=1  \text{ in } \Omega \\
   & v_0=0 \text{ on } \partial \Omega.
\end{split}
\end{equation}
To start the process,  we need that $v_0$ is smooth in the interior.  This is guaranteed by the fact that $\Omega$ is close to $B_1$.  More precisely,  we proved in \cite{CX}:
\begin{lem}\label{dirichlet problem estimate}
Let $\Omega\subset \mathbb{C}^n$ be a bounded domain and $B_{1-\gamma}(0)\subset \Omega\subset B_{1+\gamma}(0)$ for some $0<\gamma<1$.  Let $v_0$ be the solution to the Dirichlet problem in (\ref{3.1}),  then 
\begin{equation*}
|z|^2-1-3\gamma \le v_0\le |z|^2-1+3\gamma.
\end{equation*}
Moreover,  there exists $\gamma_n>0$ small enough,  such that if $\gamma\le \gamma_n$,  we have $v_0\in C^{4}(\bar{B}_{0.9})$ with $||v_0-(|z|^2-1)||_{C^4,  B_{0.9}}\le C$.  Here $C$ depends only on $n$.
\end{lem}

The following Lemma is also proved in the $W^{2,p}$ paper:

\begin{lem}\label{l5}
Assume that $\det u_{i\bar{j}}=f$ in $\Omega$ and $u|_{\partial\Omega}=0$.  Let $v_0$ be the solution to the Dirichlet problem (\ref{3.1}).  Assume that $1- \eps \le f\le 1+\eps$.  Then we have $(1+\eps)^{\frac{1}{n}}v_0\le u\le (1-\eps)^{\frac{1}{n}}v_0$.  In particular
\begin{equation*}
|v_0-u|\le 4\eps \text{ in $\Omega$.}
\end{equation*}
\end{lem}

Then we can prove the following lemma:
\begin{lem}\label{nonintersection}
There exists a constant $\delta>0$. For any $\Bar{\epsilon}\in(0,e^{-1})$, we can choose $\mu_0$ and $\epsilon$ to be small depending on $\Bar{\epsilon}$ and $n$ such that given a section $S_t(x)$ with $t\le \mu_0^2$ and $y\notin S_t(x)$, we have that
    \begin{equation*}
        B_{\epsilon_2^{\delta}}(T(y))\cap T(S_{(1-\epsilon_2)t}(x)) =\emptyset.
    \end{equation*}
    for any $\Bar{\epsilon}<\epsilon_2 <e^{-1}$. Here $T=\frac{1}{\sqrt{t}}T_{t,x}^{-1}$ is an affine transformation such that 
    \begin{equation*}
        B_{(1-\sigma)}(0) \subset \frac{1}{\sqrt{t}} T_{t,x}^{-1} (S_t(x)-x) \subset  B_{(1+\sigma)}(0).
    \end{equation*}
\end{lem}
\begin{proof}
    We prove this lemma in two cases.

    (1) $(1-\epsilon_2)t$ and $t$ are in the same generation, i.e. there exists $k$ such that $\mu_0^{k+1}< (1-\epsilon_2)t \le t \le \mu_0^k$. Define a new coordinate $w=\frac{1}{\sqrt{t}}T_{t,x}^{-1}(z-x)$. In this coordinate, $S_t(x)$ is $\sigma-$ close to the unit ball as in the Lemma \ref{section similar to a ball 2} and there is a plurisubharmonic function $h$ such that $\phi-h=0$ on $\partial S_t(x)$ by the Lemma \ref{characterization of sections}. From the argument of the section 2, we can just assume that $\phi=0$ on $\partial S_t(x)$ and use the coordinate $w$ in the rest of the proof because the linearized Monge-Ampere equation is invariant under the normalization. Using the Lemma \ref{section similar to a ball 2}, the Lemma \ref{dirichlet problem estimate} and the Lemma \ref{l5} we can get that:
    \begin{equation*}
        |w|^2-1-3\epsilon -C\sigma \le \phi \le |w|^2 -1+3\epsilon +C\sigma.
    \end{equation*}
    After the normalization, $\phi(0)=-1$. This is because the height of the section is determined by the value of $\phi$ at $0$ and the height is scaled to be 1 under the normalization (The details can be found in the construction of the sections in the $W^{2,p}$ paper). So for any $w\in \partial S_{(1-\epsilon_2)t}(x)$, we have that $\phi=-\epsilon_2$ and:
    \begin{equation*}
        |w|^2-1-3\epsilon -C\sigma \le -\epsilon_2 \le |w|^2 -1+3\epsilon +C\sigma
    \end{equation*}
    This implies that:
    \begin{equation*}
        1-3\epsilon -C\sigma -\epsilon_2 \le |w|^2 \le 1+3\epsilon +C\sigma -\epsilon_2. 
    \end{equation*}
     Since $y\notin S_t(x)$ and $(1-\sigma)B_1(0) \subset T (S_t(x)-x)$, we have that:
     \begin{equation*}
         T(y-x) \notin (1-\sigma)B_1(0).
     \end{equation*}
     This implies that
     \begin{equation*}
         B_{\epsilon_2^{\delta}}(T(y-x)) \subset B_{1-\sigma-\epsilon_2^{\delta}}^c(0),
     \end{equation*}
     Here $B_{1-\sigma-\epsilon_2^{\delta}}^c(0)=\mathbb{C}^n \setminus B_{1-\sigma-\epsilon_2^{\delta}}(0)$
     In order to make sure that
     \begin{equation*}
         B_{\epsilon_2^{\delta}}(T(y-x)) \cap B_{\sqrt{1+3\epsilon +C\sigma -\epsilon_2}}(0)= \emptyset,
     \end{equation*}
      it suffices to require that 
     \begin{equation*}
         \sqrt{1+3\epsilon +C\sigma -\epsilon_2} \le 1-\sigma-\epsilon_2^{\delta}.
     \end{equation*}
     This is equivalent to 
     \begin{equation*}
         3\epsilon +(C+2)\sigma \le \epsilon_2 +\epsilon_2^{2\delta} -2\epsilon_2^{\delta}-2\sigma \epsilon_2^{\delta}.
     \end{equation*}
     If we let $\delta$ be big(this is independent of $\Bar{\epsilon}$), then the minimum of $\epsilon_2 +\epsilon_2^{2\delta} -2\epsilon_2^{\delta}-2\sigma \epsilon_2^{\delta}$ with $\Bar{\epsilon}\le \epsilon_2 \le \frac{1}{2}$ is $\bar{\epsilon} +\bar{\epsilon}^{2\delta} -2\bar{\epsilon}^{\delta}-2\sigma \bar{\epsilon}^{\delta}$.
     So we just require that 
     \begin{equation*}
         3\epsilon +(C+2)\sigma \le \bar{\epsilon} +\bar{\epsilon}^{2\delta} -2\bar{\epsilon}^{\delta}-2\sigma \bar{\epsilon}^{\delta}.
     \end{equation*}
     Noting that $\bar{\epsilon}\le e^{-1}$. When $\delta$ is big we have that: 
     \begin{equation*}
         \bar{\epsilon} +\bar{\epsilon}^{2\delta} -2\bar{\epsilon}^{\delta}-2\sigma \bar{\epsilon}^{\delta} \ge \frac{1}{2}\bar{\epsilon}.
     \end{equation*}
     So we just require that:
     \begin{equation*}
         3 \epsilon +(C+2)\sigma \le \frac{1}{2}\bar{\epsilon}.
     \end{equation*}

     (2) $(1-\epsilon_2)t$ and $t$ are not in the same level, i.e. there exists an integer $k$ such that:
     \begin{equation*}
         \mu_0^{k+1}< (1-\epsilon_2)t \le \mu_0^{k+1} < t \le \mu_0^k.
     \end{equation*}
     Recall that the sections $S_t(x)$ for $\mu_0^{k+1}< t \le \mu_0^k$ can be written as 
     \begin{equation*}
         S_t(x)=\{z: \phi-h(z)\le \phi(x)+t\},
     \end{equation*}
     according to the Lemma \ref{characterization of sections}. We can define $\widetilde{S}_t(x)$:
     \begin{equation}\label{tilde S}
         \widetilde{S}_t(x)=\{z: \phi-h(z)\le \phi(x)+t\},
     \end{equation}
     for $t\le \mu_0^{k+1}$. Note that when we define $S_{t}(x)$ for $t\le \mu_0^{k+1}$, we subtract $\phi$ by a pluriharmonic function that is different from $h$ and then take sublevel sets. So $\widetilde{S}_t(x)$ is different from $S_t(x)$ for $t\le \mu_0^{k+1}$. Similar to the Proposition \ref{section similar to a ball 2}, we can show that 
     \begin{equation*}
         S_{(1-\epsilon_2)t}(x)\subset \widetilde{S}_{(1+c(\sigma))(1-\epsilon_2)t}(x).
     \end{equation*}
     When we let $\mu_0$ be small and $\epsilon$ be small, $\sigma$ can be arbitrarily small. Then $c(\sigma)$ can be arbitrarily small. As is defined in the case (1), we use the coordinate $w$ where $S_t(x)$ is normalized to be close to the unit ball. So using the Equation \ref{tilde S}, for any $w\in \partial \widetilde{S}_{(1+c(\sigma))(1-\epsilon_2)t}(x)$, we have that $\phi=(1+c(\sigma))(1-\epsilon_2)-1=c(\sigma)-\epsilon_2 -\epsilon_2 c(\sigma)$. Since in this coordinate
     \begin{equation*}
         |w|^2-1-3\epsilon -C\sigma \le \phi \le |w|^2-1+3\epsilon +C\sigma,
     \end{equation*}
     we have that:
     \begin{equation*}
         |w|^2-1-3\epsilon -C\sigma \le c(\sigma) -\epsilon_2 -\epsilon_1 c(\sigma) \le |w|^2 -1 +3\epsilon +C\sigma.
     \end{equation*}
     This is equivalent to that:
     \begin{equation*}
         1-3\epsilon -C\sigma +c(\sigma) -\epsilon_2 -\epsilon_2 c(\sigma) \le |w|^2 \le 1+3\epsilon +C\sigma +c(\sigma) -\epsilon_2 -\epsilon_2 c(\sigma).
     \end{equation*}
     Since $y\notin S_t(x)$ and
     \begin{equation*}
         (1-\sigma) B_1(0) \subset T( S_t(x)-x).
     \end{equation*}
     So we have that
     \begin{equation*}
         T(y-x) \notin (1-\sigma)B_1(0).
     \end{equation*}
     So we have that 
     \begin{equation*}
         B_{\epsilon_2^{\delta}}(T(y-x)) \subset B_{1-\sigma -\epsilon_2^{\delta}}^c(0).
     \end{equation*}
     So if we want 
     \begin{equation*}
         B_{\epsilon_2^{\delta}}(T(y-x)) \cap B_{\sqrt{1+3\epsilon +C\sigma -\epsilon_2 +c(\sigma) -\epsilon_2 c(\sigma)}}(0) = \emptyset,
     \end{equation*}
     we just require that
     \begin{equation*}
         B_{1-\sigma -\epsilon_2^{\delta}}^c(0)\cap B_{\sqrt{1+3\epsilon +C\sigma -\epsilon_2 +c(\sigma) -\epsilon_2 c(\sigma)}}(0) = \emptyset.
     \end{equation*}
     This is implied by
     \begin{equation*}
         \sqrt{1+3\epsilon +C\sigma -\epsilon_2 +c(\sigma) -\epsilon_2 c(\sigma)} \le 1-\sigma -\epsilon_2^{\delta}.
     \end{equation*}
     This is equivalent to 
     \begin{equation*}
         3\epsilon +C\sigma +c(\sigma) +2\sigma -\sigma^2 \le \epsilon_2 +\epsilon_2 c(\sigma) +\epsilon_2^{2 \delta}+2\epsilon_2^{\delta}+2\sigma \epsilon_2^{\delta}.
     \end{equation*}
     Let $c(\sigma)$ be small and let $\delta$ be the same as in the case (1). We can see that $f(t)=t +t c(\sigma) +t^{2 \delta}+2t^{\delta}+2\sigma t^{\delta}$ takes the minimum at $t=\Bar{\epsilon}$. So we just require that:
     \begin{equation*}
         3\epsilon +C\sigma +c(\sigma) +2\sigma -\sigma^2 \le \bar{\epsilon} +\bar{\epsilon} c(\sigma) +\bar{\epsilon}^{2 \delta}+2\bar{\epsilon}^{\delta}+2\sigma \bar{\epsilon}^{\delta}.
     \end{equation*}
     We can assume that $\delta$ is big and $\Bar{\epsilon}\in (0,\frac{1}{2})$ such that $\bar{\epsilon} +\bar{\epsilon} c(\sigma) +\bar{\epsilon}^{2 \delta}+2\bar{\epsilon}^{\delta}+2\sigma \bar{\epsilon}^{\delta}\ge \frac{1}{2}\Bar{\epsilon}$.
     So we jut require that
     \begin{equation*}
         3\epsilon +C\sigma +c(\sigma) +2\sigma -\sigma^2 \le \frac{1}{2}\Bar{\epsilon}.
     \end{equation*}
     This is true if we let $\mu_0$ be small and $\epsilon$ be small.
\end{proof}

Then we can prove the following lemma which is similar to the Lemma 1 in \cite{CG} with some modifications.

\begin{lem}\label{intersection number}
    For any $\Bar{\epsilon}<e^{-1}$, we can let $\mu_0$ be small and let $\epsilon$ be small. For any $A \subset B_{\frac{1}{2}}(0)$ which is a bounded set. Fix a positive function $\widetilde{t}$ defined on $A$ satisfying $0< \tilde{t}\le \frac{\mu_0}{2}$. Let us denote by $F$ the family of all the sections $S_{\widetilde{t}(x)}(x)$ with $x\in A$. Then there exists a countable subfamily of $F$, $\{S_{\widetilde{t}(x_k)}(x_k)\}_{k=1}^{\infty}$ with the following properties ( For simplicity we denote $\widetilde{t}(x_k)$ as $t_k$ from now on ):

    (i) For a.e. $x\in A$, $x\in \cup_{k=1}^{\infty}S_{t_k}(x_k)$.

    (ii) $x_k \notin \cup_{j<k}S_{t_j}(x_j)$, $\forall k\ge 2$.

    (iii) We have that the family 
    \begin{equation*}
        F_{\epsilon_2}=\{S_{(1-\epsilon_2)t_k}(x_k)\}_{k=1}^{\infty}
    \end{equation*}
    has bounded overlaps; more precisely
    \begin{equation*}
        \Sigma_{k=1}^{\infty} \chi_{S_{(1-\epsilon_2)t_k}(x_k)}(x)\le K log\frac{1}{\epsilon_2},
    \end{equation*}
    for $\Bar{\epsilon}\le \epsilon_2\le e^{-1}$, where $K$ is a constant which depends only on $n$.
\end{lem}
\begin{proof}
We first state the idea of the proof of the Lemma 1 in the paper \cite{CG}. Then state where we should make modifications to adapt the proof to our case.

The idea of the proof of the Lemma 1 in the paper \cite{CG} is the following: Fix a constant $M$. For each $k\ge 1$, we consider the sections $S_{t_x}(x)$ where $\frac{M}{2^{k+1}}< t_x\le \frac{M}{2^k}$. Take a countable subset of these sections, denoted as $F_k'$. 

To show (i), we shall first prove that each generation $F_k'$ has overlapping bounded by a constant that is independent of $k$. Second, we shall deduce from this that each generation $F_k'$ has a finite number of members; in particular, by relabeling the members of $F_k'$ we obtain (ii)(Here we need to use the way we define $F_k'$). This implies that the process in the construction of $F_k'$ stopped at some point and then (i) follows. For the proof of (iii), it suffices to estimate the overlapping of sections belonging to different generations. Suppose that 
\begin{equation*}
    z_0\in \cap_i S_{(1-\epsilon)t_{r_i}^{e_i}}(x_{r_i}^{e_i}),
\end{equation*}
with $S_{t_{r_i}^{e_i}}(x_{r_i}^{e_i})$ to be a section defined in the generation $i$ and $e_1<e_2<...<e_i<...$. Using the Lemma \ref{advanced engulfing property}, we can see that if $e_i-e_1$ is big $z_0$ is close to $x_{r_i}^{e_i}$. However, by the Lemma \ref{nonintersection}, the distance between $z_0$ and $x_{r_i}^{e_i}$ is bounded from below by $\epsilon_2^{\delta}$. This gives an upper bound on $e_i-e_1$.

Now we state where we should make modifications to adapt the proof to our case. Since the height of the sections that we are considering are always smaller than $1$, we can just take $M=1$ in the proof. The proof of the Lemma 1 in the paper \cite{CG} mainly use assumption (A), (B) and (C) in that paper. 

The assumption (A) corresponds to the Lemma \ref{advanced engulfing property} in our paper. Note that in the statement of Lemma \ref{advanced engulfing property}, fortunately, the balls that are used to be compared with $\frac{1}{\sqrt{t_0}}T_{t_0,x_0}^{-1}S_t(x)$ are centered exactly at $\frac{1}{\sqrt{t_0}}T_{t_0,x_0}^{-1} x$. However, the assumption (A) uses balls which may not center at $\frac{1}{\sqrt{t_0}}T_{t_0,x_0}^{-1} x$. This allows us make the proof easier when we use the Lemma \ref{advanced engulfing property}. 

The assumption (B) corresponds to the Lemma \ref{nonintersection}. In our case, we need to require additionally that $\Bar{\epsilon}<\epsilon_2<\frac{1}{2}$ to use the Lemma \ref{nonintersection} and thus we need this requirement in the statement of the Lemma we are proving. We will see in the rest of the paper that this restriction doesn't cause any trouble.

We don't have the assumption (C) in our paper because the assumption (C) is about sections with height going to the infinity while in our case the height of the sections we defined are bounded by $1$. \cite{CG} uses the assumption (C) to prove that each generation $F_k'$ is finite. Next we prove that each generation $F_k'$ is finite in our case without using the assumption (C). Denote $F_k' =\{S_{t_i^k}(x_i^k)\}$. Then we can use the Lemma \ref{diam limit} to get that:
\begin{equation*}
 B_{\frac{1}{C}(t_i^k)^{\frac{1}{2}+log_{\mu_0}(1-C\sigma^{\frac{1}{2}})}}(x_i^k) \subset S_{t_i^k}(x_i^k) \subset B_{C (t_i^k)^{\frac{1}{2}+log_{\mu_0} (1+C\sigma^{\frac{1}{2}})}}(x_i^k).
\end{equation*}
According to the proof of the Lemma 1 in the paper \cite{CG} we can prove that the intersection in each generation is bounded by a constant $\alpha^n$ that is independent of $k$ without using the assumption (C). So the intersection of $\{B_{\frac{1}{C}(t_i^k)^{\frac{1}{2}+log_{\mu_0}(1-C\sigma^{\frac{1}{2}})}}(x_i^k)\}$ is also bounded by $\alpha^n$, i.e. 
\begin{equation}\label{sigma chi estimate}
    \Sigma_i \chi_{B_{\frac{1}{C}(t_i^k)^{\frac{1}{2}+log_{\mu_0}(1-C\sigma^{\frac{1}{2}})}}(x_i^k)}(x)\le \alpha^n.
\end{equation}
Since $A\subset B_{\frac{1}{2}}(0)$, we can assume that $\mu_0$ is small and us the Lemma \ref{diam limit} to get that 
\begin{equation*}
 B_{\frac{1}{C}(t_i^k)^{\frac{1}{2}+log_{\mu_0}(1-C\sigma^{\frac{1}{2}})}}(x_i^k) \subset S_{t_i^k}(x_i^k) \subset B_{C (t_i^k)^{\frac{1}{2}+log_{\mu_0} (1+C\sigma^{\frac{1}{2}})}}(x_i^k) \subset B_1(0).
\end{equation*}
Then we integrate the Inequality \ref{sigma chi estimate} on $B_1(0)$ to get that:
\begin{equation*}
\begin{split}
    &\Sigma_i \alpha(2n)(\frac{1}{C}(t_i^k)^{\frac{1}{2}+log_{\mu_0}(1-C\sigma^{\frac{1}{2}})})^{2n} \\
    &=\Sigma_i \int_{B_1(0)}  \Sigma_i \chi_{B_{\frac{1}{C}(t_i^k)^{\frac{1}{2}+log_{\mu_0}(1-C\sigma^{\frac{1}{2}})}}(x_i^k)}(x) dx \le \alpha(2n) \alpha^n,
\end{split}
\end{equation*}
where $\alpha(2n)$ is the volume of the unit ball in $\mathbb{C}^n$. Note that $\frac{1}{2^{k+1}}<t_i^k \le \frac{1}{2^k}$ and we assume that $0< log_{\mu_0}(1-C\sigma^{\frac{1}{2}})<\frac{1}{2}$. we can get that 
\begin{equation*}
    \Sigma_i \alpha(2n)(\frac{1}{C}(\frac{1}{2^{k+1}}))^{2n}\le \alpha(2n)\alpha^n.
\end{equation*}
This implies that the number of terms in the first sum is finite. So there is only finite sections in $F_k'$.
\end{proof}

Next we want to estimate the difference between $\mu(S_t(x))$ and $\mu(S_{(1-\epsilon_2)t}(x))$.
\begin{lem}\label{smaller sections error}
    For $t\le \mu_0$, there exists a constant $C$ depending only on $n$ such that 
    \begin{equation*}
        \mu(S_t(x) \setminus S_{(1-\epsilon_2)t}(x)) +\mu(S_{(1-\epsilon_2)t}(x) \setminus S_t(x))\le C(\sigma +\epsilon_2)\mu(S_t(x)).
    \end{equation*}
    Moreover, we have that:
    \begin{equation*}
        |\mu(S_t(x))-\mu(S_{(1-\epsilon_2)t}(x))|\le C(\sigma +\epsilon_2)\mu(S_t(x)).
    \end{equation*}
\end{lem}
\begin{proof}
    Using the Lemma \ref{section similar to a ball 2}, we have that 
    \begin{equation*}
        \begin{split}
            &B_{(1-\sigma)}(0)\subset \frac{1}{\sqrt{t}}T_{t,x}^{-1}S_t(x) \subset B_{(1+\sigma)}(0) \\
            & B_{(1-\sigma)}(0)\subset \frac{1}{\sqrt{(1-\epsilon_2)t}}T_{(1-\epsilon_2)t,x}^{-1}S_{(1-\epsilon_2)t}(x) \subset B_{(1+\sigma)}(0)
        \end{split}
    \end{equation*}
    Since for any $t\le \mu_0$, the determinant of $T_{t,x}$ is one, so it is area preserving. Then we have that 
    \begin{equation*}
    \begin{split}
               & \alpha(2n)(1-\sigma)^{2n}t^n \le m(S_t(x)) \le \alpha(2n)(1+\sigma)^{2n}t^n  \\
               & \alpha(2n)(1-\sigma)^{2n}((1-\epsilon_2)t)^n \le m(S_{(1-\epsilon_2)t}(x)) \le \alpha(2n)(1+\sigma)^{2n}((1-\epsilon_2)t)^n.
    \end{split}
    \end{equation*}
    So we have that 
    \begin{equation*}
        m(S_t(x) \setminus S_{(1-\epsilon_2)t}(x)) +m(S_{(1-\epsilon_2)t}(x) \setminus S_t(x))\le C(\sigma +\epsilon_2)m(S_t(x)).
    \end{equation*}
    SInce $|f-1|\le \epsilon$, we have that:
    \begin{equation*}
        \frac{1}{2}m \le (1-\epsilon)m \le \mu \le (1+\epsilon)m \le 2m.
    \end{equation*}
    So we have that 
    \begin{equation*}
        \mu(S_t(x) \setminus S_{(1-\epsilon_2)t}(x)) +\mu(S_{(1-\epsilon_2)t}(x) \setminus S_t(x))\le C(\sigma +\epsilon_2)\mu(S_t(x)).
    \end{equation*}
    
\end{proof}

\begin{proof}
(of Theorem \ref{cz})
   The proof of this theorem is similar to the Theorem in \cite{CG}. We will sketch the proof and point out the modifications we need to make for our case. 
   
   The idea of the proof is the following: Since $\frac{\mu(S_k \cap A)}{\mu(S_k)}\le \delta$, intuitively we have that: 
   \begin{equation*}
       \mu(A)\le \delta \mu(\cup_{k=1}^{\infty} S_k).
   \end{equation*}
   This is true if $S_k$ don't intersect with each other. Although in general this can't happen, we can consider smaller sections: $S_k^{\epsilon_2}=S_{(1-\epsilon_2)t_k}(x_k)$. The overlapping of the smaller sections can be estimated by the Lemma \ref{nonintersection}. Note that the norm between $S_k$ and $S_k^{\epsilon_2}$ is quite small by the Lemma \ref{smaller sections error}. So the error caused by using $S_k^{\epsilon_2}$ instead of $S_k$ is under control.

   Now we want to point out the modifications we need to make for our case. 

The Theorem in the \cite{CG} only requires that the sections $S_k$ satisfy a weaker property:
\begin{equation*}
    \frac{\delta}{C_1}\le \frac{\mu(S_k\cap A)}{\mu(S)}\le \delta.
\end{equation*}
instead of the property (a) that we want to prove. We indeed need the stronger version in the rest of the paper. We need to define $S_k$ in a more careful way. For each $x\in A$ such that $t_x$ is well defined as in the Lemma \ref{intersection number}, we consider $\nu_x=\lambda_x t_x$ where $1< \lambda_x<2$ is to be determined. Using the definition of $t_x$ and the left continuity of the section with respect to its height, we have that:
\begin{equation*}
    \frac{\mu(S_{t_x}(x)\cap A)}{\mu(S_{t_x}(x))}\ge \delta.
\end{equation*}
We also have that:
\begin{equation*}
    \frac{\mu(S_t(x)\cap A)}{\mu(S_t(x))}\le \delta, \text{ for any } t_x\le t \le \frac{t_x}{\mu_0}.
\end{equation*}
This is because by definition, for any $t_x\le t\le \mu_0^3$, 
\begin{equation*}
    \frac{\mu(S_t(x)\cap A)}{\mu(S_t(x))}\ge \delta
\end{equation*}
and $t_x\le \mu_0^4$.
Then we have that:
\begin{equation}\label{mu ratio}
    \begin{split}
        & \delta \ge \frac{\mu(S_{t}(x)\cap A)}{\mu(S_{t}(x))}\ge \frac{\mu(S_{t_x}(x)\cap A)}{\mu(S_{t}(x))} \\
        &\ge \frac{\mu(S_{t_x}(x)\cap A)}{\mu(S_{t_x}(x))}(1-C\sigma-C(\frac{t}{t_x}-1)) \ge (1-C\sigma-C(\frac{t}{t_x}-1)) \delta
    \end{split}
\end{equation}
for any $(1+c(\sigma))t_x\le (1+C\sigma^{\frac{1}{2}})t_x< t\le 2t_x$ if we assume that $\mu_0\le \frac{1}{2}$. Note that we use $(1+c(\sigma))t_x\le (1+C\sigma^{\frac{1}{2}})t_x< t$ and the Proposition \ref{section similar to a ball 2} to get $S_{t_x}(x) \subset S_t(x)$. In the second line we use the Lemma \ref{smaller sections error}. As in the proof of the Theorem in \cite{CG}, we want to apply the Lemma \ref{intersection number} to $\{S_{\nu_x}(x)\}$ to get a subfamily of the sections, denoted as $\{S_k=S_{\mu_k}(x_k)\}$. We also define $S_k^{\epsilon_2}=S_{(1-\epsilon_2)\mu_k}(x_k)$. Then go through the calculation of the proof of the Theorem in \cite{CG}. There is one thing that requires modification in the calculation. \cite{CG} uses the following inequality:
 \begin{equation*}
     |\mu(S_{\nu_k}(x_k))-\mu(S_{(1-\epsilon_2)\nu_k}(x_k))|\le \epsilon_2 logC \mu(S_{\nu_k}(x_k)).
 \end{equation*}
 In our case, we can use the Lemma \ref{smaller sections error} to prove that:
    \begin{equation*}
        |\mu(S_{\nu_k}(x_k))-\mu(S_{(1-\epsilon_2)\nu_k}(x_k))|\le C(\sigma +\epsilon_2)\mu(S_{\nu_k}(x_k)).
    \end{equation*}
 This is a little weaker than \cite{CG}. However, we will see it is enough for the rest of the proof. In our case, we can eventually get that:
 \begin{equation*}
     \mu(A)=\mu((\cup S_k)\cap A)\le \theta(\epsilon_2)\mu(\cup S_k^{\epsilon_2}),
 \end{equation*}
   where 
   \begin{equation*}
       \theta(\epsilon_2)=(1-\frac{1-\delta}{C_0 log(\frac{1}{\epsilon_2})})(1+\frac{C(\sigma +\epsilon_2)}{1-C(\gamma+\epsilon_2)}C_0 log(\frac{1}{\epsilon_2})).
   \end{equation*}
The constant $C$ in the above formula depends only on $n$. As in the \cite{CG}, for any $0<\epsilon_2<C(n,\delta)$ we have that:
 \begin{equation*}
     (1-\frac{1-\delta}{C_0 log(\frac{1}{\epsilon_2})})(1+\frac{2C\epsilon_2}{1-2C\epsilon_2}C_0 log(\frac{1}{\epsilon_2}))<1
 \end{equation*}
Then we fix $\epsilon_2\in (0,C(n,\delta))$ and let $\mu_0$ be small and $\epsilon$ be small such that $\sigma\le \epsilon_2$ and the $\Bar{\epsilon}$ in the Lemma \ref{intersection number} is smaller than $\epsilon_2$. So we have that:
\begin{equation*}
    \theta(\epsilon_2)<1.
\end{equation*}
The property (c) is satisfied. We want $S_k$ to satisfy the property (a). We use the Inequalities \ref{mu ratio} to get that:
\begin{equation*}
        \delta \ge \frac{\mu(S_{\lambda_x t_x}(x)\cap A)}{\mu(S_{\lambda_x t_x}(x))}\ge (1-C\sigma-C(\lambda_x-1)) \delta
\end{equation*}
Recall that we need $\lambda_x\ge 1+C\sigma^{\frac{1}{2}}$ to get the Inequalities \ref{mu ratio}. So we let $\lambda_x =1+C\sigma^{\frac{1}{2}}$ so that:
\begin{equation*}
    (1-C\sigma-C(\lambda_x-1)) \delta \ge (1-C\sigma^{\frac{1}{2}} )\delta.
\end{equation*}
The property (a) is satisfied.
\end{proof}

\section{Critical density}
The Theorem 1 in \cite{CG2} also holds in our case:
\begin{thm}\label{critical density}
    There are two constants $M_1>1$ and $0<\lambda <1$, depending only on the doubling constant of $\mu$ and dimension, such that for any section $S=S_t(x)$ and any nonnegative solution $u$ of $L_{\varphi}u=0$ such that
    \begin{equation*}
        \inf_{z\in S_{\frac{t}{2}}(x)} u(z)\le 1
    \end{equation*}
    we have that
    \begin{equation*}
        \mu (\{z\in S: u(z)>M_1\}) < \lambda \mu(S).
    \end{equation*}
\end{thm}
\begin{proof}
    The proof of the theorem is similar to the Theorem 1 in \cite{CG2}. We first sketch the proof. First we normalize $S_t(x)$ to be close to the unit ball and normalize $\phi$ correspondingly which will be made clearer latter. Let $\alpha$ be big enough such that $v=u+\alpha \phi$ satisfies that 
    \begin{equation*}
        \inf_{z\in S_{\frac{t}{2}}(x)}v(z)\le -1.
    \end{equation*}
    We will show that we can take $\alpha=6$ in our case. Then we can use the ABP estimate for $v^-(x)=-\min{\{v(z),0\}}$ to get an lower bound estimate on the measure of the points where $v^-$ and its convex envelope $\Gamma(v^-)$ touch. On these touching points, $u$ is uniformly bounded because at such point, $v\le 0$ and $\phi$ is uniformly bounded by the ABP estimate and $\alpha$ is a uniform constant which we will explain later. This concludes the proof of the theorem.

    Next we point out the modifications we use to adapt the proof of the Theorem 1 in \cite{CG2} to our case. 
    
    First we normalize $\phi$ and change the coordinate (denote the new coordinate as $w$). According to the section 2, the Proposition \ref{section similar to a ball 2}, the Lemma \ref{characterization of sections}, we can assume that:
    \begin{equation*}
        \phi=0 \text{ on } \partial S_t(x)
    \end{equation*}
    and 
    \begin{equation*}
        (1-0.1\sigma)B_1(0) \subset S_t(x) \subset (1+0.1\sigma)B_1(0).
    \end{equation*}
     Claim:  We can let $\sigma$ and $\epsilon$ be small enough such that:
    \begin{equation*}
        \phi\le -\frac{1}{3} \text{ on } S_{\frac{t}{2}}(x).
    \end{equation*}
    We prove this in two cases. Before the proof, we note that by the Proposition \ref{section similar to a ball 2}, we have that there exists an ellipsoid $E$ such that:
    \begin{equation*}
        (1-0.1\sigma)\sqrt{\frac{1}{2}}E \subset S_{\frac{t}{2}}(x) \subset (1+0.1\sigma)\sqrt{\frac{1}{2}}E.
    \end{equation*}

    Case (1) $S_t(x)$ and $S_{\frac{t}{2}}(x)$ are in the same generation. So we have that $E=B_1(0)$. So we have that:
    \begin{equation*}
        B_{\frac{1-0.1\sigma}{\sqrt{2}}}(0) \subset S_{\frac{t}{2}}(x) \subset B_{\frac{1+0.1\sigma}{\sqrt{2}}}(0).
    \end{equation*}
    Using the Lemma \ref{section similar to a ball 2}, the Lemma \ref{dirichlet problem estimate} and the Lemma \ref{l5} we can get that:
    \begin{equation}\label{phi estimate}
        |w|^2-1-3\epsilon -C\sigma \le \phi \le |w|^2 -1+3\epsilon +C\sigma.
    \end{equation} 
    So we have that for any $w\in S_{\frac{t}{2}}(x)$, 
    \begin{equation*}
        \phi\le  \frac{(1+0.1\sigma)^2}{2}-1+3\epsilon +C\sigma.
    \end{equation*}
    Let $\sigma$ and $\epsilon$ be small enough, we can get that 
    \begin{equation*}
        \phi\le -\frac{1}{3} \text{ on } S_{\frac{t}{2}}(x).
    \end{equation*}

    Case (2) $S_t(x)$ and $S_{\frac{t}{2}}(x)$ are not in the same generation. We can assume that $\mu_0<\frac{1}{2}$ so that there exists $k$ such that:
    \begin{equation*}
        \mu_0^{k+2} <\frac{t}{2}\le \mu_0^{k+1} < t \le \mu_0^k,
    \end{equation*}
    i.e. $S_t(x)$ and $S_{\frac{t}{2}}(x)$ are in adjacent generations. So we have that:
    \begin{equation*}
        B_{(1-C\sigma^{\frac{1}{2}})}(0) \subset E \subset B_{(1+C\sigma^{\frac{1}{2}})}(0),
    \end{equation*}
    according to the Proposition \ref{section similar to a ball 2} and the Lemma \ref{tx0k estimate}. So we have that:
    \begin{equation*}
        B_{\frac{(1-0.1\sigma)(1-C\sigma^{\frac{1}{2}})}{\sqrt{2}}}(0) \subset S_{\frac{t}{2}}(x) \subset B_{\frac{(1+0.1\sigma)(1+C\sigma^{\frac{1}{2}})}{\sqrt{2}}}(0).
    \end{equation*}
    Using the Inequalities \ref{phi estimate}, we have that for any $w\in S_{\frac{t}{2}}(x)$,
    \begin{equation*}
        \phi\le \frac{(1+0.1\sigma)^2(1+C\sigma^{\frac{1}{2}})^2}{2}-1+3\epsilon +C\sigma.
    \end{equation*}
     Let $\sigma$ and $\epsilon$ be small enough, we can get that 
    \begin{equation*}
        \phi\le -\frac{1}{3} \text{ on } S_{\frac{t}{2}}(x).
    \end{equation*}
    The claim is proved.
    
    The ABP estimate needs to use $det(D^2 \phi)$. In our case we only know the estimate for $det (\phi_{i\Bar{j}})$ instead of $det(D^2 \phi)$. However, we can use $det(D^2 \phi)^{\frac{1}{2}}\le 2^n det(\phi_{i\bar{j}})$ at the points where $\phi$ is convex. Fortunately the integral part in the ABP estimate is taken only on the points where $\phi$ is convex.

    The constant $\alpha$ we use is uniform. This is because after normalization, we have that $\phi=0$ on $\partial S_t(x)$ and $\phi\le -\frac{1}{3}$ on $S_{\frac{t}{2}}(x)$. In the assumption we have that $\inf_{z\in S_{\frac{t}{2}}(x)}u(z)\le 1$. So we can just take $\alpha=6$ so that 
    \begin{equation*}
        \inf_{z\in S_{\frac{t}{2}}(x)}(u(z)+\alpha \phi(z))\le -1.
    \end{equation*}
\end{proof}

\section{Infimum estimate on a larger section}
We need the following Lemma:
\begin{lem}\label{Luepsilon}
    Let $u$ be a positive supersolution of the equation $L_{\phi}u=0$, i.e.,
    \begin{equation*}
        L_{\phi}u \le 0.
    \end{equation*}
    If $x_0\in R^n$ and $0<\epsilon_4 <1$ then 
    \begin{equation*}
        L_{\phi}(-u^{\epsilon_4})(x_0) \ge -\epsilon_4(\epsilon_4-1)u(x_0)^{\epsilon_4-2}|\nabla u (x_0)|^2 \frac{det (u_{i\Bar{j}}(x_0))}{trace(u_{i\Bar{j}})}.
    \end{equation*}
\end{lem}
\begin{proof}
    The lemma can be proved by following the proof of the Lemma 2.1 in \cite{CG2} word by word.
\end{proof}

\begin{thm}\label{larger set inf estimate}
There is a constant $L>1$ depending on $n$, such that for any $\alpha>0$, if $u$ is any nonnegative solution of $L_{\phi}u =0$ on $S_{4t}(x_0)$ such that:
\begin{equation*}
    u|_{S_t(x_0)}> \alpha,
\end{equation*}
then
\begin{equation*}
    u|_{S_{2 t}(x_0)}>\frac{\alpha}{L}
\end{equation*}
\end{thm}
\begin{proof}
This theorem is similar to the Theorem 2 in the paper \cite{CG2}. We sketch the proof and point out the modifications we need to adapt the proof of the Theorem 2 in the paper \cite{CG2} to our case.

First we sketch the proof of the theorem as follows: Consider four sections: $S_k^*=S_{t_k}(x_0)$, $k=1,2,3,4$, where $t_1=t<t_2=2t <t_3<t_4=4t$. $t_3$ is to be determined (This is the definition in the paper \cite{CG2}. We will change the definition of $S_3^*$ in our case.). We normalize $\phi$ and change the coordinate (denote the new coordinate as $w$) such that:
\begin{equation*}
    \phi|_{\partial S_{4t}(x_0)}=0
\end{equation*}
and
\begin{equation*}
    B_{1-0.1\sigma}(0) \subset S_{4t}(x_0) \subset B_{1+0.1\sigma}(0).
\end{equation*}
Then define an auxiliary function $w_{\epsilon_4}$ whose definition will be talked about in detail latter. Then consider 
\begin{equation*}
    h(x)=\frac{-\phi(x)+w_{\epsilon_4}(x)}{2\beta_n},
\end{equation*}
where $\beta_n$ is a constant depending on $n$. Consider the minimum of $u^{\epsilon_4} -h(x)$ on $S_4^*$ and prove that
\begin{equation*}
    u|_{S_{2t}(x)}\ge \frac{\alpha}{L}
\end{equation*}
in three cases:

Case 1. The minimum is attained on $S_1^*$.

Case 2. The minimum is attained on $S_4^*\setminus S_3^*$.

Case 3. The minimum is attained on $S_3^* \setminus S_1^*$.

The case 1 and 2 don't even need the auxiliary function $w_{\epsilon_4}$. The case 3 is more intricate. One can estimate $L_{\phi}(u^{\epsilon_4} -h(x))$ at the minimum point using the Lemma \ref{Luepsilon}. Since $L_{\phi}(u^{\epsilon_4} -h(x))\ge 0$ at the minimum point, we can show that at this point $g=\Delta \phi$ is bounded which can be used to get further estimate.

Next we talk about the modifications that we need to adapt the proof of the Theorem 2 in the paper \cite{CG2}. 

First, $S_3^*$ in our case may not take the form of $S_t(x)$. Instead, we define 
\begin{equation*}
    S_3^*=\{w\in S_4^*:\phi(w)\le s_3\},
\end{equation*}
where $s_3\in[-\frac{3}{8},-\frac{1}{4}]$ is to be determined. Using the Lemma \ref{section similar to a ball 2}, the Lemma \ref{dirichlet problem estimate} and the Lemma \ref{l5} we can get that:
    \begin{equation}\label{phi estimate}
        |w|^2-1-3\epsilon -C\sigma \le \phi \le |w|^2 -1+3\epsilon +C\sigma.
    \end{equation} 
So we have that for any $\phi \le s_3$:
\begin{equation*}
|w|^2\le 1+3\epsilon +C\sigma +s_3.
\end{equation*}
Then we have that
\begin{equation*}
    S_3^* \subset B_{\sqrt{1+3\epsilon+C\sigma +s_3}}(0) \subset B_{\sqrt{1+3\epsilon+C\sigma -\frac{1}{4}}}(0) \subset B_{\frac{2}{\sqrt{5}}}(0).
\end{equation*}
Here we assume that $\sigma$ and $\epsilon$ are small.

In \cite{CG2} an estimate for $|\nabla \phi|$ in $S_3^*$ is needed. We want to derive such estimate in our case. We have shown that 
\begin{equation*}
    S_3^* \subset B_{\frac{2.1}{\sqrt{5}}}(0).
\end{equation*}
In our case we can use the result in the $W^{2,p}$ estimate: $||\phi||_{W^{2,p}(B_{\frac{2}{\sqrt{5}}}(0))}\le C$. If $p$ is big then we have that $||\phi||_{C^{1,\beta_0}(B_{\frac{2}{\sqrt{5}}}(0))}\le C$ for some $\beta_0>0$. In particular, we get an estimate for $|\nabla \phi|$ in $S_3^*$.

In \cite{CG2}, an estimate of the measure of the set is needed: 
\begin{equation*}
    H_{\epsilon_4}=\{x\in S_3^*: g(x)\ge \frac{\gamma_0}{\epsilon_4}\},
\end{equation*}
where $\gamma_0$ is a positive small parameter. The estimate of the area of $\partial S_3^*$ is needed which is easy when $S_3^*$ is convex as in \cite{CG2}. In our case $S_3^*$ is just a pseudoconvex set. So we use the coarea formula:
\begin{equation*}
    \int_{\{a\le \phi \le b\}} |\nabla \phi| dx =\int_a^b dt\int_{\phi=t}dS,
\end{equation*}
for any $a<b$. By the sard's theorem, we can find $s_{ab}\in (a,b)$ such that $\{\phi=s_{ab}\}$ is smooth and
\begin{equation*}
    \int_a^b dt\int_{\phi=t}dS \ge (b-a)|\{\phi=s_{ab}\}|.
\end{equation*}
We need to estimate 
\begin{equation*}
    m(\{a \le \phi\le b\}).
\end{equation*}

 Using the Lemma \ref{section similar to a ball 2}, the Lemma \ref{dirichlet problem estimate} and the Lemma \ref{l5} we can get that:
    \begin{equation}\label{phi estimate}
        |w|^2-1-3\epsilon -C\sigma \le \phi \le |w|^2 -1+3\epsilon +C\sigma.
    \end{equation} 
So we have that for any $a\le \phi \le b$:
\begin{equation*}
    1-3\epsilon-C\sigma +a \le |w|^2\le 1+3\epsilon +C\sigma +b.
\end{equation*}

So we have that:
\begin{equation*}
    \{a\le \phi \le b\}\subset B_{\sqrt{1+3\epsilon+C\sigma+b}}(0) \setminus B_{\sqrt{1-3\epsilon-C\sigma+a}}(0).
\end{equation*}
We assume that $b\le -\frac{1}{4}$. We have proved that $|\nabla \phi|_{\{\phi\le -\frac{1}{4}\}}\le C$, we can combining the above equations to get that
\begin{equation}\label{level set estimate}
\begin{split}
        &(b-a)|\{\phi=s_{ab}\}|\le \int_{\{a\le \phi \le b\}} |\nabla \phi| dx  \\
        &\le |\nabla \phi|_{C^1(\{a\le \phi\le b\})} m(B_{\sqrt{1+3\epsilon+C\sigma+b}}(0) \setminus B_{\sqrt{1-3\epsilon-C\sigma+a}}(0)).
\end{split}
\end{equation}

We use the Lemma \ref{smaller sections error} to get that:
\begin{equation*}
    m(B_{\sqrt{1+3\epsilon+C\sigma+b}}(0) \setminus B_{\sqrt{1-3\epsilon-C\sigma+a}}(0))\le C (\sigma+\epsilon +b-a).
\end{equation*}
We can use this inequality and the Inequality \ref{level set estimate} to get that:
\begin{equation*}
    |\{\phi=s_{ab}\}|\le C(1+\frac{\sigma+\epsilon}{b-a}),
\end{equation*}
for some $s_{ab}\in [a,b]$. So when we want to select $t_3$, we can use the above argument to find $s_3\in[-\frac{3}{8},-\frac{1}{4}]$. Then we define  From the calculation above, we have that
\begin{equation*}
    |\{\phi=s_3\}|\le C(1+\sigma+\epsilon).
\end{equation*}
Using this we can estimate $|H_{\epsilon_4}|$ as in \cite{CG2} to get:
\begin{equation*}
    |H_{\epsilon_4}|\le A_n \frac{\epsilon_4}{\gamma_0} (1+\sigma+\epsilon) \le C\frac{\epsilon_4}{\gamma_0}
\end{equation*}

Another modification that we want to make is the following: In \cite{CG2}, an auxiliary function is defined. We first state the way \cite{CG2} define an auxiliary function: Let
\begin{equation*}
    k(x)=det D^2\phi(x)
\end{equation*}
Then approximate the set $H_{\epsilon_4}$ by an open set $\tilde{H}_{\epsilon_4}$ such that $H_{\epsilon_4}\subset \tilde{H}_{\epsilon_4} \subset S_4^*$ such that $|\tilde{H}_{\epsilon_4}\setminus H_{\epsilon_4}|$ is sufficiently small. Given $\delta>0$(small), they define $\varphi(x)$ a smooth function in $S_4^*$ such that $\varphi(x)=1$ in $H_{\epsilon_4}$, $\varphi(x)=\delta $ in $S_4^* \setminus \tilde{H}_{\epsilon_4}$, and $\delta \le \varphi(x) \le 1$. Then they define the auxiliary function:
\begin{equation*}
    \begin{split}
        &det D^2 w_{\epsilon_4}(x) =k(x) \varphi(x) \text{ in } S_4^* \\
        & w_{\epsilon_4}|_{\partial S_4^*}=0.
    \end{split}
\end{equation*}

Next we explain how to define an auxiliary function in our case. Since we are dealing with complex Monge-Ampere equation and linearzed complex Monge-Ampere equation, it is natural to define the auxiliary function using complex Monge-Ampere equation. However, later on we need to derive an estimate for the derivative of the auxiliary function. Such estimate is missing for the complex Monge-Ampere equation. So we still want to use the real Monge-Ampere equation to define the auxiliary function. Another difficulty is that in our case $S_4^*$ may not be convex. So we can't solve the dirichlet problem for the real Monge-Ampere equation on $S_4^*$. So instead, we solve the dirichlet problem on an ellipsoid which is slightly larger than $S_4^*$:
\begin{equation*}
    \begin{split}
        &det D^2 w_{\epsilon_4}(x) =4^n k^2(x) \varphi^2(x) \text{ in } (1+\sigma)B_1(0) \\
        & w_{\epsilon_4}|_{(1+\sigma)B_1(0)}=0,
    \end{split}
\end{equation*}
where $k(x)=det \phi_{i\Bar{j}}$. Recall that from the Proposition \ref{section similar to a ball 2},
\begin{equation*}
    (1-\sigma)B_1(0) \subset S_4^* \subset (1+\sigma)B_1(0).
\end{equation*}
 If we let $\epsilon_4$ be small and let $\delta$ be small, we can get that $|w_{\epsilon_4}|$ and $|D(w_{\epsilon_4})|_{S_3^*}$ are very small as in \cite{CG2}. Then we consider the three cases as listed above.
 
(1) We can prove the case 1 word by word as in \cite{CG2}.

(2)For the second case. Since $S_{4(1+s)t}(x_0)$ and $S_{4t}(x_0)$ are in the same or the adjacent generations for $-\frac{1}{2}\le s\le 0$, we have that in the coordinate $w$:
\begin{equation}\label{S B compare}
    (1-C\sigma^{\frac{1}{2}})B_{\sqrt{1+s}}(0) \subset S_{4(1+s)t}(x_0) \subset (1+C\sigma^{\frac{1}{2}})B_{\sqrt{1+s}}(0).
\end{equation}
Combining this with the Inequalities \ref{phi estimate}, we have that:
\begin{equation*}
    \{\phi\le -\frac{1}{2}-C(\sigma^{\frac{1}{2}}+\epsilon)\}\subset S_2^* \subset \{\phi\le -\frac{1}{2}+C(\sigma^{\frac{1}{2}}+\epsilon)\}.
\end{equation*}
By the definition of $s_3$,
\begin{equation*}
    \{\phi\le -\frac{3}{8}\} \subset S_3^* \subset \{\phi\le -\frac{1}{4}\}
\end{equation*}
 So for $x\in S_2^*$ we have
\begin{equation*}
\begin{split}
        &h(x)-h(P)=\frac{-\phi^*(x)}{2\beta_n} +\frac{w_{\epsilon_4}(x)}{2\beta_n} +\frac{\phi^*(P)}{2\beta_n} -\frac{w_{\epsilon_4}(P)}{2\beta_n}\\
        &\ge \frac{\frac{1}{2}-C(\sigma^{\frac{1}{2}}+\epsilon)}{2\beta_n} +\frac{w_{\epsilon_4}(x)}{2\beta_n} +\frac{-\frac{3}{8}}{2\beta_n} \\
        &=\frac{1-C(\sigma^{\frac{1}{2}}+\epsilon)}{16\beta_n} +\frac{w_{\epsilon_4}(x)}{2\beta_n} \\
        &\ge \frac{1}{32\beta_n}
\end{split}
\end{equation*}
In the last line we let $\sigma$, $\epsilon$ and $\epsilon_4$ be small enough.

(3) For the case three.
As in \cite{CG2}, we need to show:
\begin{equation*}
    |\nabla \phi^* (P)|\ge C>0,
\end{equation*}
for $P\in S_3^*- S_1^*$.  Recall that 
\begin{equation*}
    S_3^* \subset B_{\frac{2}{\sqrt{5}}}(0).
\end{equation*}
We need to use a version of interpolation (see the Lemma 4.32 of \cite{GT}). We can use that 
\begin{equation*}
    |\phi-(|z|^2-1)|_{B_{\sqrt{\frac{5}{6}}}(0)}\le C(\epsilon+\sigma).
\end{equation*}
and 
\begin{equation*}
|\phi|_{C^{1,\alpha}(B_{\sqrt{\frac{5}{6}}}(0))} \le C    
\end{equation*}
by the $W^{2,p}$ estimate. Then using interpolation and get that $|\nabla (\phi -(|z|^2-1))|_{B_{\frac{2}{\sqrt{5}}}(0)}$ can be arbitrarily small if $\sigma$ and $\epsilon$ are small enough. In particular, we can assume that:
\begin{equation*}
    |\nabla (\phi -(|z|^2-1))|_{B_{\frac{2}{\sqrt{5}}}(0)} \le \frac{1}{\sqrt{3}}.
\end{equation*}
Using the Formulae \ref{S B compare}, we have that 
\begin{equation*}
    S_2^* \supset B_{\frac{1}{\sqrt{3}}}(0).
\end{equation*}
Note that:
\begin{equation*}
    |\nabla (|z|^2-1)|_{B_{\frac{1}{\sqrt{3}}}^c (0)}\ge \frac{2}{\sqrt{3}}.
\end{equation*}
So we have that:
\begin{equation*}
    |\nabla \phi|_{S_3^* \setminus S_2^*} \ge \frac{1}{\sqrt{3}}.
\end{equation*}
This can imply
\begin{equation*}
    |\nabla \phi^* (P)|\ge \frac{1}{\sqrt{3}}.
\end{equation*}

We also need to use the inequality:
\begin{equation*}
    det D^2 w_{\epsilon_4}\le 4^n det ((w_{\epsilon_4})_{i\Bar{j}})^2.
\end{equation*}
to reduce $det ((w_{\epsilon_4})_{i\Bar{j}})$ to $det D^2 w_{\epsilon_4}$. Those are all the modifications that we need.
\end{proof}

\section{Harnack inequality and H\"older continuity}
First we can combine the Theorem \ref{critical density} and the Theorem \ref{larger set inf estimate} to get the following Lemma:
\begin{lem}\label{stronger critical density}
    Let $u$ be a nonnegative solution in a section $S_{4t}(x_0)$, and $\alpha$ is any positive number such that:
    \begin{equation*}
        \mu\{x\in S_t(x_0):u(x)>\alpha\}\ge \lambda \mu(S_t(x_0)),
    \end{equation*}
    then 
    \begin{equation*}
        u(x) \ge \frac{\alpha}{M_0}, \,\,\, \forall x\in S_t(x_0),
    \end{equation*}
    Here $M_0=M_1 L^2$. $\lambda$ and $M_1$ are given in the Theorem \ref{critical density} and $L$ is given in the Theorem \ref{larger set inf estimate}.
\end{lem}
\begin{proof}
We can define $v=\frac{M_1 u}{\alpha}$. $v$ still satisfies the linearized comlex Monge-Ampere equation. Using the Theorem \ref{critical density}, we have that
\begin{equation*}
\inf_{x\in S_{\frac{t}{2}}(x_0)}v(x)\ge 1. 
\end{equation*}
Since $2^2\ge 2(1+c(\sigma))$. Using the Proposition \ref{section similar to a ball 2}, we have that
\begin{equation*}
    S_{2t}(x_0)\supset S_t (x_0).
\end{equation*}
Using the Therorem \ref{larger set inf estimate}, we have that:
\begin{equation*}
    \inf_{S_{2t}(x_0)}v \ge \frac{1}{L}\inf_{S_{t}(x_0)}v \ge \frac{1}{L^2}\inf_{S_{\frac{1}{2}t}(x_0)}v\ge \frac{1}{L^2}
\end{equation*}
Then we have that:
\begin{equation*}
    \inf_{S_t(x_0)}u \ge \inf_{S_{2t}(x_0)}u =\frac{\alpha}{M_1}\inf_{S_{2t}(x_0)}v \ge \frac{\alpha}{M_1 L^2}
\end{equation*}
\end{proof}

The Lemma 3.1 of \cite{CG2} also holds in our case with $t \le \mu_0^3$:
\begin{lem}\label{height estimate}
    Let $u$ be a nonnegative solution of $L_{\phi}u=0$ in the section $S_{t}(z)$ such that:
    \begin{equation*}
        \inf_{S_{t}(z)}u \le 1
    \end{equation*}
    with $t\le \mu_0^3$. Let $\theta$ be the constant in the Lemma \ref{engulfing property theta version}. Then if $y\in S_{t}(z)$ and $S_h(y)$ is a section with $h< \theta t$, and
    \begin{equation*}
        \mu\{x\in S_h(y):u(x)>\alpha\}\ge \lambda \mu(S_h(y))
    \end{equation*}
    then 
    \begin{equation*}
        h\le \theta t (\frac{M_0 L}{\alpha})^{\frac{1}{\delta}}.
    \end{equation*}
    Here, $0<\lambda$ is given in the Theorem \ref{critical density}. $M_0$ is given by the Lemma \ref{stronger critical density}. $\theta$ is the constant in the Lemma \ref{engulfing property theta version}, and $L>1$ is the constant in the Theorem \ref{larger set inf estimate}.
\end{lem}
\begin{proof}
    We can prove this Lemma following the proof of the Lemma 3.1 of \cite{CG2} word by word.
\end{proof}

The following Lemma comes from \cite{CX}. 
\begin{lem}\label{L0.10}
Let $f:B_{0.8}\rightarrow \mathbb{R}$ be an $L^1$ function.  Then for $m$-a.e.  $x\in B_{0.8}$,  we have:
\begin{equation*}
\lim\sup_{x\in S_{\mu_{\alpha}}(x_{\alpha}),\,\mu_{\alpha}\rightarrow 0}\frac{1}{m(S_{\mu_{\alpha}}(x_{\alpha}))}\int_{S_{\mu_{\alpha}}(x_{\alpha})}|f(y)-f(x)|dm(y)=0.
\end{equation*}
In particular, for any $A\subset B_{0.8}$, we can take $f=\chi_A$ in the above formula and get that: For a.e. $x\in A$,
\begin{equation*}
    \lim \sup_{x\in S_{\mu_{\alpha}}(x_{\alpha}),\,\mu_{\alpha}\rightarrow 0}|\frac{1}{m(S_{\mu_{\alpha}}(x_{\alpha}))}m(A\cap S_{\mu_{\alpha}}(x_{\alpha}))-1|=0.
\end{equation*}
\end{lem}

Next we want to prove a decay estimate for the suplevel sets of $u$ and then get $L^p$ estimate for $u$ for some $p>0$.
\begin{thm}\label{lp inf}
    There exists $p>0$ such that for any $z_0\in B_{\frac{1}{2}}$ we have that:
    \begin{equation*}
        |u|_{L^p(B_{\frac{\mu_0^2}{C}}(z_0))}\le C \inf_{S_{\mu_0^3}(z_0)} u.
    \end{equation*}
    Where $B_{\frac{\mu_0^2}{C}}(z_0)\subset S_{\mu_0^3}(z_0)$.
\end{thm}
\begin{proof} We will do the calculation in a set $S_{\mu_0^3}(z_0)$. By multiplying $u$ by a constant which doesn't affect the conclusion, We may assume that 
\begin{equation*}
\inf_{S_{\mu_0^3}(z_0)}u \le 1.
\end{equation*}
Let $B_{\frac{\mu_0^2}{C}}(z_0)$ be a subset of $S_{\mu_0^3}(z_0)$. This is ensured by the Lemma \ref{diam limit}. We define 
\begin{equation*}
    E_k=\{x: u(x)\ge KM^k\},
\end{equation*}
where $K$ and $M$ are constants to be determined. Define
\begin{equation*}
    S_k=B_{\frac{\mu_0^2}{C}(1-\frac{1}{4}-...-\frac{1}{2^{k+1}})}(z_0)
\end{equation*}
for $k\ge 1$ and $S_0=B_{\frac{\mu_0^2}{C}}(z_0)$.
We want to prove that:
\begin{equation*}
    \mu(E_{k+1}\cap S_k) \le \delta_1 \mu(E_k \cap S_{k-1})
\end{equation*}
for some $0<\delta_1<1$.
Before we get the Calderon-Zygmund decomposition of $S_1\cap E_2$, we need to verify that for a.e. $x\in A$, the assumptions (1) and (2) in the Theorem \ref{cz} hold. The assumption (1) holds because of the Lemma \ref{L0.10}. Suppose that the assumption (2) is not true. Then there exists $t_0\in (\mu_0^4, \mu_0^3]$ such that:
\begin{equation*}
    \mu(S_{t_0}(x)\cap E_2 \cap S_1)>\lambda \mu(S_{t_0}(x)).
\end{equation*}
Then we can use the Lemma \ref{stronger critical density} to get that:
\begin{equation*}
    \inf_{x\in S_{t_0}(x)} u(x)\ge \frac{KM^2}{M_0}.
\end{equation*}

Since
\begin{equation*}
\inf_{S_{\mu_0^3}(z_0)}u \le 1
\end{equation*}
and $t_0\le \mu_0^3$, we can use the Lemma \ref{height estimate} to get that
\begin{equation*}
    t_0\le \theta \mu_0^3 (\frac{M_0 L}{KM^2})^{\frac{1}{\delta}}.
\end{equation*}
Let $M$ be big such that $t_0\le \mu_0^4$. This is a contradiction because we already assume that $t_0> \mu_0^4$. So the assumptions (1) and (2) in the Theorem \ref{cz} hold.

Then we can use the Theorem \ref{cz} to get a Calderon-Zygmund decomposition of $S_1\cap E_2$ at the level $\frac{\lambda}{1-C\sigma^{\frac{1}{2}}}$: $\{S_{t_i^{1}}(x_i^{(1)})\}$ (We assume that $\sigma$ and $\epsilon_2$ are small such that $\frac{\lambda}{1-C\sigma^{\frac{1}{2}}}<1$).
Moreover, we have that:
\begin{equation*}
    \frac{\mu(S_{t_i^{1}}(x_i^{(1)}))\cap S_1 \cap E_2}{\mu(S_{t_i^{1}}(x_i^{(1)}))}\ge (1-C\sigma^{\frac{1}{2}})\frac{\lambda}{1-C\sigma^{\frac{1}{2}}} =\lambda.
\end{equation*}
Then we can use the Theorem \ref{stronger critical density} to get that:
\begin{equation*}
    u\ge \frac{KM^2}{M_0}\ge KM \text{ on } S_{t_i^{(1)}}(x_i^{(1)}).
\end{equation*}
Here we assume that $M\ge M_0$. This implies that
\begin{equation*}
    S_{t_i^{(1)}}(x_i^{(1)}) \subset E_1.
\end{equation*}
Using the Lemma \ref{height estimate}, we can get that:
\begin{equation*}
    t_i^{(1)}\le \theta \mu_0^3(\frac{M_0 L}{KM^{2}})^{\frac{1}{\delta}}.
\end{equation*}
This implies that
\begin{equation*}
    S_{t_i^{(1)}}(x_i^{(1)}) \subset B_{[\theta \mu_0^3(\frac{M_0 L}{KM^{2}})^{\frac{1}{\delta}}]^{\frac{1}{2}+log_{\mu_0} (1+C\sigma^{\frac{1}{2}})}}(x_i^{(1)}),
\end{equation*}
by the Lemma \ref{diam limit}. We can let $M$ be big enough such that:
\begin{equation*}
    [\theta \mu_0^3(\frac{M_0 L}{KM^{2}})^{\frac{1}{\delta}}]^{\frac{1}{2}+log_{\mu_0} (1+C\sigma^{\frac{1}{2}})}\le \frac{1}{2^2}\frac{\mu_0^2}{C}.
\end{equation*}
This implies that:
\begin{equation*}
    S_{t_i^{(1)}}(x_i^{(1)}) \subset S_0.
\end{equation*}
So we have proved that 
\begin{equation*}
    S_{t_i^{(1)}}(x_i^{(1)}) \subset S_0 \cap E_1.
\end{equation*}
Using the Theorem \ref{cz}, we can get that 
\begin{equation*}
    \mu(S_1 \cap E_2) \le \delta_0 \mu(\cup  S_{t_i^{(1)}}(x_i^{(1)})) \le \delta_0 \mu(S_0 \cap E_1).
\end{equation*}
Suppose that we have proved that:
\begin{equation*}
    \mu(S_{k-1}\cap E_{k}) \le \delta_0 \mu(S_{k-2}\cap E_{k-1}).
\end{equation*}
Follow the proof above we can get a Calder\'on-Zygmund decomposition for $S_k \cap E_{k+1}$ at the level $\frac{\lambda}{1-C\sigma^{\frac{1}{2}}}$ from the Theorem \ref{cz}: $\{S_{t_i^{(k)}}(x_i^{(k)})\}$. We have that:
\begin{equation}\label{decay}
    \mu(S_{k}\cap E_{k+1}) \le \delta_0 \mu(\cup S_{t_i^{(k)}}(x_i^{(k)}))
\end{equation}
Moreover, the Theorem \ref{cz} implies that
\begin{equation*}
    \frac{\mu(S_{t_i^{(k)}}(x_i^{(k)}) \cap E_{k+1})}{\mu(S_{t_i^{(k)}}(x_i^{(k)}))}\ge (1-C\sigma^{\frac{1}{2}})\frac{\lambda}{1-C\sigma^{\frac{1}{2}}} =\lambda.
\end{equation*}
Then we can use the Lemma \ref{stronger critical density} to get that:
\begin{equation*}
    u\ge \frac{KM^{k+1}}{M_0}\ge KM^k \text{ on } S_{t_i^{(k)}}(x_i^{(k)}).
\end{equation*}
This implies that
\begin{equation*}
    S_{t_i^{(k)}}(x_i^{(k)}) \subset E_k.
\end{equation*}
Using the Theorem \ref{height estimate}, we can get that:
\begin{equation*}
    t_i^{(k)}\le \theta \mu_0^3(\frac{M_0 L}{KM^{k+1}})^{\frac{1}{\delta}}.
\end{equation*}
This implies that
\begin{equation*}
    S_{t_i^{(k)}}(x_i^{(k)}) \subset B_{[\theta \mu_0^3(\frac{M_0 L}{KM^{k+1}})^{\frac{1}{\delta}}]^{\frac{1}{2}+log_{\mu_0}(1+C\sigma^{\frac{1}{2}})}}(x_i^{(k)}),
\end{equation*}
by using the diameter estimate.
We can let $M$ be big enough which is independent of $k$ such that:
\begin{equation*}
    [\theta \mu_0^3(\frac{M_0 L}{KM^{k+1}})^{\frac{1}{\delta}}]^{\frac{1}{2}+log_{\mu_0}(1+C\sigma^{\frac{1}{2}})}\le \frac{1}{2^{k+1}}\frac{\mu_0^2}{C}.
\end{equation*}
This implies that:
\begin{equation*}
    S_{t_i^{(k)}}(x_i^{(k)}) \subset S_{k-1}.
\end{equation*}
So we have proved that 
\begin{equation*}
    S_{t_i^{(k)}}(x_i^{(k)}) \subset S_{k-1} \cap E_k.
\end{equation*}
So the Inequality \ref{decay} implies that:
\begin{equation*}
    \mu(S_{k} \cap E_{k+1})\le \delta_0 \mu(S_{k-1} \cap E_k).
\end{equation*}
So we have that:
\begin{equation*}
    \mu(B_{\frac{\mu_0^2}{2C}}(z_0)\cap E_{k+1})\le \mu(S_k \cap E_{k+1})\le \delta_0^{k} \mu(S_0 \cap E_1).
\end{equation*}
So we have that:
\begin{equation*}
\begin{split}
        &\int_{B_{\frac{\mu_0^2}{2C}}(z_0)}u^p =p\int_0^{\infty} s^{p-1}Area(\{x:u(x)\ge s\}\cap B_{\frac{\mu_0^2}{2C}}(z_0)) ds \\
        &=p\int_0^{KM}s^{p-1}ds +\Sigma_{i=1}^{\infty} \int_{KM^i}^{KM^{i+1}}s^{p-1}Area(\{x:u(x)\ge s\}\cap B_{\frac{\mu_0}{2C}}(z_0))ds \\
        & \le C+\Sigma_{i=1}^{\infty} p\int_{K M^{i}}^{KM^{i+1}} s^{p-1} Area(\{u\ge KM^{i}\}\cap B_{\frac{\mu_0^2}{2C}}(z_0))ds \\
        &\le C+ C\Sigma_{i=1}^{\infty}(KM^{i+1})^p \delta_0^{i-1}.
\end{split}
\end{equation*}
Then we can take $p$ be small such that $M^p \delta_0 <1$. So the last line in the above formula is finite. 
\end{proof}

Next we need a Lemma which is proved in \cite{CC}:
\begin{lem}\label{sup lp}
    Suppose that $u\ge 0$ satisfies in $B_1\subset R^d$:
    \begin{equation*}
        \partial_i (a^{ij}\partial_j u)\ge 0.
    \end{equation*}
    Here $\frac{1}{\lambda(x)}\le a^{ij}(x)\le \lambda(x)$, with $\lambda(x)\in L^p(B_1)$ for some $p>\frac{3d}{2}$, then for any $\epsilon_5>0$, there exists a constant $C$, depending on $p$, $||\lambda||_{L^p(B_1)}$ and $\epsilon_5$ such that:
    \begin{equation*}
        \sup_{B_{\frac{1}{2}}}u \le C||u||_{L^{\epsilon_5}(B_1)}.
    \end{equation*}
\end{lem}

Combining the Theorem \ref{lp inf} and the Lemma \ref{sup lp}, we can prove the following Harnack inequality:
\begin{thm}\label{initial harnack}
    There exists a constant $C$ such that for any $z_0\in B_{\frac{1}{2}}$, we have that:
    \begin{equation*}
        \sup_{B_{\frac{\mu_0^2}{C}}(z_0)}u \le C\inf_{B_{\frac{\mu_0^2}{C}}(z_0)}u
    \end{equation*}
\end{thm}
\begin{proof}
    Using the main theorem of the $w^{2,p}$ paper, we can derive the $W^{2,p}$ estimate for $\phi$ and $p$ can be made arbitrarily big if we let $\epsilon$ be small enough. This implies that $\max_l{\lambda_l(\phi_{i\Bar{j}})}\in L^p(B_{0.8})$. Since $1-\epsilon\le det (\phi_{i\Bar{j}})\le 1+\epsilon$, we can have that 
    \begin{equation*}
        \max_k{\frac{1}{\lambda_k(\phi_{i\Bar{j}})}} \le 2\max_l{\lambda_l(\phi_{i\Bar{j}})}^{n-1}
    \end{equation*}
    Thus we can also get that $\max_k{\frac{1}{\lambda_k(\phi_{i\Bar{j}})}} \in L^p(B_{0.8})$. So we can use the Lemma \ref{sup lp} to get that:
    \begin{equation*}
       \sup_{B_{\frac{\mu_0^2}{2C}}(z_0)}u \le C||u||_{L^{\epsilon_5}(B_{\frac{\mu_0^2}{C}}(z_0))}.
    \end{equation*}
    Using the Lemma \ref{lp inf}, we have that:
      \begin{equation*}
        ||u||_{L^{p_0}(B_{\frac{\mu_0^2}{C}}(z_0))}\le C \inf_{S_{\mu_0^3}(z_0)} u \le C\inf_{B_{\frac{\mu_0^2}{2C}}(z_0)}u
    \end{equation*}
    for some $p_0>0$. Then we can let $\epsilon_5$ be small enough such that $\epsilon_5\le p_0$. Then we have that:
    \begin{equation*}
        \sup_{B_{\frac{\mu_0^2}{2C}}(z_0)}u \le C||u||_{L^{\epsilon_5}(B_{\frac{\mu_0^2}{C}}(z_0))} \le C ||u||_{L^{p_0}(B_{\frac{\mu_0^2}{C}}(z_0))} \le C\inf_{B_{\frac{\mu_0^2}{2C}}(z_0)}u.
    \end{equation*}
\end{proof}

We should be able to prove the Theorem \ref{harnack all scales}:

\begin{proof}
(of the Theorem \ref{harnack all scales}.)
    We first normalize $S_{t}(z_0)$ to be close to a ball. Define a coordinate $w$ by $z-z_0=\sqrt{t}T_{t,z_0}w$. Let $h_{z_0}$ be a degree two pluriharmonic polynomial such that when we define 
    \begin{equation*}
    \begin{split}
        &\widetilde{\phi}(w)=\frac{\phi(z_0+\sqrt{t}T_{t,z_0} w)}{t |det(T_{t,z_0})|^{\frac{2}{n}}} +h_{z_0},\\
        &\widetilde{u}(w)=u(z_0+\sqrt{t}T_{t,z_0} w)
        \end{split}
    \end{equation*}
    we have that $\widetilde{\phi}=0$ on $\partial S_{t}(z_0)$ in the coordinate $w$. Using the Theorem 1.3 in the $W^{2,p}$ paper, we can get that 
    \begin{equation*}
||\widetilde{\phi}||_{W^{2,p}(B_{0.9})}\le C,  \,\,\,\,\,||\sum_i\frac{1}{\widetilde{\phi}_{i\bar{i}}}||_{L^p(B_{0.9})}\le C,
\end{equation*}
for any big $p$ if we let $\mu_0$ be small and $\epsilon$ be small. Then we can apply the Theorem \ref{initial harnack} to get that:
\begin{equation*}
    \sup_{B_{\frac{\mu_0^2}{C}}} \widetilde{u} \le C\inf_{B_{\frac{\mu_0^2}{C}}} \widetilde{u}
\end{equation*}
Let $\widetilde{S}_{\mu}(w_0)$ be the sections defined in the coordinate $w$ in the same way that we define $S_{\mu}$. From the Lemma \ref{section similar to a ball 2}, we have that:
\begin{equation*}
     \frac{1}{1-c_1(\gamma)}\widetilde{S}_{\frac{\mu_0^4 t}{C_0}}(0) \subset B_{\frac{\mu_0^2}{C}},
\end{equation*}
if we let $C_0$ be big enough and assume that $c_1(\gamma)\le \frac{1}{2}$. Using the above formula and the relation between $u$ and $\widetilde{u}$, we have that:
\begin{equation*}
    \sup_{\frac{1}{1-c_1(\gamma)}\widetilde{S}_{\frac{\mu_0^4 t}{C_0}}(0) } u \le C\inf_{\frac{1}{1-c_1(\gamma)}\widetilde{S}_{\frac{\mu_0^4 t}{C_0}}(0) } u.
\end{equation*}
Then use the Lemma \ref{key lemma 2} to get that:
\begin{equation*}
    S_{\frac{\mu_0^4 t}{C_0}}(z_0) \subset \frac{1}{1-c_1(\gamma)}\widetilde{S}_{\frac{\mu_0^4 t}{C_0}}(0).
\end{equation*}
So we have proved that:
\begin{equation*}
    \sup_{S_{\frac{\mu_0^4 t}{C_0}}(z_0) } u \le C\inf_{S_{\frac{\mu_0^4t}{C_0}}(z_0)} u.
\end{equation*}
\end{proof}

\begin{proof}
(of the Corollary \ref{main theorem baby})
 Define $\tau=\frac{\mu_0^4}{C_0}$. Define $\Bar{u}=u-\inf u$ so that $\Bar{u}$ is nonnegative and we can apply the Harnack inequality to $\Bar{u}$. Define $M_t(z)=\sup_{S_t(z)}\Bar{u}$ and $m_t (z) =\inf_{S_t(z)}\Bar{u}$. We can use the Corollary \ref{harnack all scales} with $u$ replaced by $M_t(z)-\Bar{u}$ and $\Bar{u}-m_t(z)$ for any $t\le \frac{\mu_0^5}{C_0}$ to get that:
    \begin{equation*}
        \begin{split}
            &M_t(z)-m_{\tau t}(z) \le \beta(M_t(z)-M_{\tau t}(z))\\
            & M_{\tau t}(z)-m_t(z)\le \beta(m_{\tau t}(z)-m_t(z)).
        \end{split}
    \end{equation*}
    Add these two inequalities together, we have that:
    \begin{equation*}
        M_{\tau t}(z)-m_{\tau t}(z) \le \frac{\beta -1}{\beta+1} (M_t (z) -m_t (z)).
    \end{equation*}
    This implies that 
    \begin{equation*}
        osc_{S_{\frac{\tau^k\mu_0^5}{C_0}}(z)}\Bar{u}\le (\frac{\beta-1}{\beta+1})^k osc_{S_{\frac{\mu_0^5}{C_0}}(z)}\Bar{u} \le 2 (\frac{\beta-1}{\beta+1})^k |\Bar{u}|_{L^{\infty}}.
    \end{equation*}
    Using the Lemma \ref{diam limit}, we have that 
    \begin{equation*}
        B_{\frac{1}{C} (\frac{\tau^k \mu_0^5}{C_0})^{\frac{1}{2}+log_{\mu_0}(1-C\sigma^{\frac{1}{2}})}}(z) \subset S_{\frac{\tau^k \mu_0^5}{C_0}}(z).
    \end{equation*}
    So we have that:
    \begin{equation*}
        osc_{B_{C (\frac{\tau^k \mu_0^5}{C_0})^{1+log_{\mu_0}(1-C\sigma^{\frac{1}{2}})}}(z)}\Bar{u} \le 2 (\frac{\beta-1}{\beta+1})^k |\Bar{u}|_{L^{\infty}}\le 4(\frac{\beta-1}{\beta+1})^k |u|_{L^{\infty}}.
    \end{equation*}
    Since $\frac{\beta-1}{\beta+1}<1$, this implies that:
    \begin{equation*}
        osc_{B_r(z)}u =osc_{B_r(z)}\Bar{u} \le C r^{\alpha},
    \end{equation*}
    which concludes the proof of the lemma.
\end{proof}

\section{Some corollaries of the main theorem}
\begin{proof}
    (of the Corollary \ref{c1.2})

    Now we take some point $p_0\in M$ and take normal coordinates $(z_1,\cdots,  z_n)$ at $p_0$ so that $g_{i\bar{j}}(p_0)=\delta_{ij}$ and $\nabla g(p_0)=0$.  We can choose local potential $\rho(z)$,  such that $\omega_0=\sqrt{-1}\partial\bar{\partial}\rho$ near $p_0$,  say on $B_1(p_0)$ (under local coordinates $z$).  So that on this neighborhood,  the equation can be written as:
\begin{equation}\label{7.1}
\det((\rho+\phi)_{i\bar{j}})=
f\det(g_{i\bar{j}}),\,\,\,\text{ in $B_1$.}
\end{equation}
In order to use Theorem \ref{main theorem},  we need to zoom in (\ref{7.1}) at $p_0$ at a suitable scale so that the right hand side is close to a contant.  

Let $0<r_0<1$,  we perform a change of variable $z=r_0w$.  Next we define 
\begin{equation*}
\tilde{u}_{r_0}(w)=\frac{1}{r_0^2}u(r_0w),\,\,\,\tilde{\rho}_{r_0}=\frac{1}{r_0^2}\rho(r_0w),\,\,\tilde{\phi}_{r_0}(w)=\frac{1}{r_0^2}\phi(r_0w).
\end{equation*}
So we have that :
\begin{equation*}
    L_{\Tilde{\rho}_{r_0}+\tilde{\phi}_{r_0}}(\tilde{u}_{r_0}-\tilde{\phi}_{r_0})=0.
\end{equation*}
From the proof of the Corollary 1.1 in the $W^{2,p}$ paper, we have that if we let $r_0$ be small enough depending on $\omega_0$ and $\epsilon$ be small, The assumptions of the Theorem \ref{main theorem} hold with $\phi$ replaced by $\Tilde{\rho}_{r_0}+\tilde{\phi}_{r_0}$ and with $u$ replaced by $\tilde{u}_{r_0}-\tilde{\phi}_{r_0}$. Then we can use the Theorem \ref{main theorem} to get that:
\begin{equation*}
    ||\tilde{u}_{r_0}-\tilde{\phi}_{r_0}||_{C^{\alpha}(B_{\frac{1}{2}}(0))}\le C.
\end{equation*}
This implies that:
\begin{equation*}
    ||\tilde{u}_{r_0}||_{C^{\alpha}(B_{\frac{1}{2}}(0))}\le C.
\end{equation*}
Using an elementary covering argument, we have that:
\begin{equation*}
    ||u||_{C^{\alpha}(M)}\le C.
\end{equation*}
\end{proof}

\end{document}